\documentclass{article}

\usepackage{amsmath}
\usepackage{amsfonts}
\usepackage{epsfig}
\usepackage{theorem}
\usepackage{comment}
\usepackage[english]{babel}
\usepackage{euscript}

\newtheorem{theorem}{Theorem}[section]
\newtheorem{lemma}[theorem]{Lemma}

\newtheorem{proposition}[theorem]{Proposition}

{\theorembodyfont{\rmfamily} 

\newtheorem{remark}[theorem]{Remark}
\newtheorem{example}[theorem]{Example}

\newtheorem{definition}[theorem]{Definition}

\newtheorem{assumption}[theorem]{Assumption}

}

\def\C{{\mathbb C}}

\def\P{{\mathbb P}}
\def\Z{{\mathbb Z}}
\def\R{{\mathbb R}}
\def\Q{{\mathbb Q}}

\newcommand\al{\alpha}

\newcommand\ga{\gamma}

\newcommand\be{\beta}
\def\eps{\varepsilon}
\newcommand\Cbar{\overline{\C}}
\newcommand\Qbar{\overline{\Q}}

\newcommand\eq{\,=\,}

\def\bsk{\bigskip}
\def\msk{\medskip}
\def\ssk{\smallskip}
\def\ni{\noindent}

\def\wt{\widetilde}

\DeclareMathOperator*{\rav}{=}

\title{Minimum Degree of the Difference of \\ Two Polynomials over~$\Q$. \\
Part~II: Davenport--Zannier pairs}

\author{Fedor Pakovich and Alexander K. Zvonkin}

\date{\today}

\begin{document}

\maketitle

\begin{abstract}
In this paper we study pairs of polynomials with a given factorization pattern and 
such that the degree of their difference attains its minimum. We call such pairs of 
polynomials {\em Davenport--Zannier pairs}, or DZ-pairs for short.
The paper is devoted to the study of DZ-pairs {\em with rational coefficients}. 

In our earlier paper \cite{PakZvo.I-14}, in the framework of the {\em theory of dessins
d'enfants}, we established a correspondence between 
DZ-pairs and {\em weighted bicolored plane trees}. These are bicolored plane trees
whose edges are endowed with positive integral weights. When such a tree is uniquely
determined by the set of black and white degrees of its vertices, it is called
{\em unitree}, and the corresponding DZ-pair is defined over~$\Q$. In~\cite{PakZvo.I-14},
we classified all unitrees. In this paper, we compute all the corresponding
polynomials. In this way, the present paper is a sequel of~\cite{PakZvo.I-14}.

In the final part of the paper we present some additional material concerning the 
Galois theory of DZ-pairs and weighted trees. 
\end{abstract}


\section{Introduction}

Let $\al,\be\vdash n$ be two partitions of an integer~$n$, 
$$
\al \eq (\al_1,\ldots,\al_p), \qquad \be \eq (\be_1,\ldots,\be_q), \qquad 
\sum_{i=1}^p\al_i \eq \sum_{j=1}^q\be_j \eq n,
$$
and let $P$ and $Q$ be two {\em coprime}\/ polynomials of degree $n$ having the following 
factorization patterns:
\begin{eqnarray}\label{eq:P-and-Q}
P(x) \eq \prod_{i=1}^p \, (x-a_i)^{\al_i}, \qquad
Q(x) \eq \prod_{j=1}^q \, (x-b_j)^{\be_j}.
\end{eqnarray}
In these expressions we consider the multiplicities $\al_i$ and $\be_j$, 
$i=1,2,\ldots,p$, $j=1,2,\ldots,q$ as being given, while the roots 
$a_i$ and $b_j$ are not fixed, though they must all be distinct. In this
paper we study polynomials satisfying (\ref{eq:P-and-Q}) and such that
{\em the degree of their difference $R = P-Q$ attains its minimum}.
Numerous papers, mainly in number theory, were devoted to the study of such
polynomials.

\begin{assumption}[Conditions on $\al$ and $\be$]\label{assump}
Throughout the paper, we always assume that
\begin{itemize}
\item   the greatest common divisor of the numbers $\al_1,\ldots,\al_p,\be_1,\ldots,\be_q$
        is~1;
\item   $p+q\le n+1$.
\end{itemize}
\end{assumption}
The case of partitions $\al$, $\be$ not satisfying the above conditions can easily be
reduced to this case (see \cite{PakZvo.I-14}).

In 1995, Zannier \cite{Zannier-95} proved that under the above conditions the following
statements hold:
\begin{enumerate}
\item   $\deg R \ge (n+1)-(p+q)$.
\item   This bound is always attained, whatever are $\al$ and $\be$.
\end{enumerate}

\begin{definition}[DZ-pair and its passport]\label{def:DZ}
A\, pair\, of\, polynomials\, $(P,Q)$ such that $P$ and $Q$ are of the form~\eqref{eq:P-and-Q} 
and $\deg\,(P-Q) = (n+1)-(p+q)$ is called {\em Davenport--Zannier pair}, or {\em DZ-pair}\/ 
for short. The pair of partitions $(\al,\be)$ is called the {\em passport}\/
of the DZ-pair.
\end{definition}

Obviously, if $(P,Q)$ is a DZ-pair with a passport $(\al,\be)$, and if we take 
$\wt{P}=c\cdot P(ax+b)$, $\wt{Q}=c\cdot Q(ax+b)$ where $ac\ne 0$, then 
$(\wt{P},\wt{Q})$ is also a DZ-pair with the same passport. We call such DZ-pairs 
equivalent.

\begin{definition}[Defined over $\Q$]\label{def:over}
We say that a DZ-pair $(P,Q)$ is {\em defined over}\/ $\Q$ if $P,Q\in\Q[x]$.
We say that an equivalence class of DZ-pairs is defined over $\Q$ if there exists 
a representative of this class which is defined over~$\Q$.
\end{definition}

By abuse of language, in what follows, we will use the shorter term ``DZ-pair'' 
to denote also an equivalence class of DZ-pairs.

\msk

In our previous paper \cite{PakZvo.I-14}, using the {\em theory of dessins d'enfants}\/ 
(see, for example, Ch.~2 of~\cite{LanZvo-04}), we established a correspondence between
DZ-pairs and {\em weighted bicolored plane trees}. These are bicolored plane trees
whose edges are endowed with positive integral weights. The degree of a vertex is
defined as the sum of the weights of the edges incident to this vertex. Obviously,
the sum of the degrees of black vertices and the sum of the degrees of white vertices
are both equal to the total weight of the tree. Let $\al=(\al_1,\al_2,\ldots,\al_p)$
and $\be=(\be_1,\be_2,\ldots,\be_q)$ be two partitions of the total weight $n$
which represent the degrees of black and white vertices respectively. The pair
$(\al,\be)$ is called the {\em passport}\/ of the tree in question.

\begin{proposition}[DZ-pairs and weighted trees]\label{biject}
There is a bijection between DZ-pairs with a passport $(\al,\be)$ on one hand, 
and weighted bicolored plane trees with the same passport on the other hand. 
\end{proposition}

\begin{definition}[Unitree]\label{uni}
A weighted bicolored plane tree such that there is no other tree with the
same passport is called {\em unitree}.
\end{definition}


General facts of the theory of dessins d'enfants imply that DZ-pairs
corresponding to unitrees are defined over $\Q$. Basing on our experience,
we claim that this class represents a vast majority of DZ-pairs defined over~$\Q$.
The other examples may roughly be subdivided into two categories. The members
of the first one are constructed as compositions of DZ-pairs corresponding 
to unitrees. The second category is, in a way, a collection of exceptions.
Still, the latter category is no less interesting since it involves some subtle 
combinatorial and group-theoretic invariants of the Galois action on DZ-pairs 
and on weighted trees.

\msk

{\bf The main result of \cite{PakZvo.I-14} is the classification of all unitrees. 
The main result of the present paper is a complete list of the corresponding
polynomials.} The final part of \cite{PakZvo.I-14} is devoted to the study of
Galois invariants of weighed trees. In the final part of the present paper we
compute the corresponding polynomials.

\msk

The class of unitrees comprises ten infinite series, denoted from $A$ to $J$, 
and ten sporadic trees, denoted from $K$ to $T$. The pictures of these trees 
are given below in the text. DZ-pairs corresponding to the series from $A$ to $J$ 
are presented in Sects.~\ref{sec:A} to \ref{sec:J}; those corresponding to the 
sporadic trees from $K$ to~$T$, in Sect.~\ref{sec:sporadic}. The Galois action
is treated in Sects.~\ref{sec:galois} to \ref{sec:other}.

For individual DZ-pairs, a computation may turn out to be difficult, sometimes 
even extremely difficult, but the verification of the result is completely trivial.
As to the infinite series, the difficulties grow as a snowball.
The ``computational'' part now consists in finding an analytic expression
of the polynomials in question, depending on one or several parameters, while the
``verification'' part consists in a {\em proof}, which may be rather elaborate.
See a more detailed discussion below.

\section{Preliminaries}

\subsection{A brief history of the question}

In 1965, Birch, Chowla, Hall, and Schinzel \cite{BCHS-65} asked a question 
which soon became famous:\label{BCHS}
\begin{quote}
Let $A$ and $B$ be two coprime polynomials with complex coefficients; 
what is the possible minimum degree of the difference $R=A^3-B^2$\,? 
\end{quote}
In order for the question to be meaningful we should take $A^3$ and $B^2$ of the same 
degree and with the same leading coefficient. Denote $\deg A = 2k$, $\deg B = 3k$, 
so that $\deg A^3 = \deg B^2 = 6k$. Let us start with an example.

\begin{example}\label{ex:pair-T}
In this example, $k=4$, so that both polynomials $P$ and $Q$ are of degree $6k=24$. 
As to their difference $R=P-Q$, all its coefficients of degrees from 24 down to~6 vanish,
so that $R$ becomes a polynomial of degree~5.
\begin{eqnarray}
P & = & (x^8+84x^6+176x^5+2366x^4+13\,536x^3+26\,884x^2 \nonumber \\
  &   &   +\,\,218\,864x+268\,777)^3, \label{P-T} \\
Q & = & (x^{12}+126x^{10}+264x^9+6195x^8+31\,392x^7+163\,956x^6 \nonumber \\
  &   &   +\,\,1\,260\,528x^5+3\,531\,639x^4+19\,770\,400x^3, \nonumber \\
  &   &   +\,\,62\,912\,622x^2+94\,024\,776x+291\,742\,453)^2, \label{Q-T} \\
R & = & -2^{38}\cdot 3^3\,(x^5+62x^3+148x^2+1001x+8852). \label{R-T}
\end{eqnarray}
\end{example}

The following two conjectures were proposed in~\cite{BCHS-65}:\label{init-problem}

\begin{enumerate}
\item   For $\deg A=2k$, $\deg B=3k$, one always has\, $\deg (A^3-B^2) \ge k+1$. 
\item   This bound is sharp: that is, it is attained for infinitely many
        values of $k$. 
\end{enumerate}

The first conjecture was proved the same year by Davenport 
\cite{Davenport-65}. The second one turned out to be much more difficult
and remained open for 16 years: in 1981 Stothers \cite{Stothers-81}
showed that the bound is in fact attained not only for infinitely many
values of $k$ but for all of them. 

A far-reaching generalization of the above result was proved in 1995
by Zannier~\cite{Zannier-95}. Let $\al=(\al_1,\ldots,\al_p)$ and
$\be=(\be_1,\ldots,\be_q)$ be two partitions of an integer $n$ satisfying
the conditions of Assumption~\ref{assump}, and let $P$ and $Q$ be two
polynomials of degree~$n$ having the factorization pattern \eqref{eq:P-and-Q}.
Then
\begin{enumerate}
\item   $\deg (P-Q) \ge (n+1)-(p+q)$.
\item   This bound is always attained, whatever are $\al$ and $\be$.
\end{enumerate}
For the case of cubes and squares considered above we have $n=6k$, 
$$
\al \eq (\underbrace{3,3,\ldots,3}_{2k}) \eq 3^{2k}, \qquad
\be \eq (\underbrace{2,2,\ldots,2}_{3k}) \eq 2^{3k},
$$ 
so that $p=2k$ and $q=3k$, whence 
$$
(n+1)-(p+q) \eq (6k+1)-(2k+3k) \eq k+1.
$$

A result equivalent to that of Zannier was, in fact, proved, in a very implicit way,
by Boccara in 1982 \cite{Boccara-82} (see also \cite{EdKuSt-84}, page 775). 
The result of \cite{Boccara-82} was purely combinatorial, and relations between 
combinatorics and polynomials were at the time largely overlooked.

Recall that a pair of polynomials $(P,Q)$ satisfying \eqref{eq:P-and-Q} and such
that the degree of $P-Q$ is {\em equal}\/ to the minimum value $(n+1)-(p+q)$ are 
called Davenport--Zannier pairs or DZ-pairs (Definition~\ref{def:DZ}). The theory
of dessins d'enfants implies that DZ-pairs are always defined over the field $\Qbar$ 
of algebraic numbers. However, the most interesting case is, without doubt, the one 
of pairs defined over $\Q$.
In 2010, Beukers and Stewart~\cite{BeuSte-10} undertook a study of DZ-pairs
of the special type $P=A^s$, $Q=B^t$, defined over~$\Q$. In our paper we study 
DZ-pairs of a general form \eqref{eq:P-and-Q} defined over $\Q$.


\subsection{Dessins d'enfants}

As we have already said,
the framework of our paper is the theory of dessins d'enfants (see, for 
example, Ch.~2 of~\cite{LanZvo-04}). The main notion of this theory is that
of Belyi function. For a rational function $f:\Cbar\to\Cbar:x\mapsto y$, where 
$\Cbar=\C\cup\{\infty\}$ is the Riemann complex sphere, let us call 
$y\in\Cbar$ a {\em critical value}\/ of~$f$ if the equation $f(x)=y$ 
has multiple roots. The definition of a Belyi function restricted to the 
planar case is as follows:

\begin{definition}[Belyi function]\label{def:belyi}
A rational function $f:\Cbar\to\Cbar$ is a {\em Belyi function}\/ if $f$ has at most
three critical values, namely, 0, 1 and $\infty$.
\end{definition}


\begin{theorem}[Belyi functions and maps]\label{th:Riemann}
If $f:\Cbar\to\Cbar:x\mapsto y$ is a Belyi function then:
\begin{itemize}
\item[\rm 1.]	The preimage ${\cal M}=f^{-1}([0,1])$ is a plane map, that is,
		a connected graph, which is embedded into the sphere in such a way that 
		its edges do not intersect. 
\item[\rm 2.]	The map ${\cal M}$ has a natural bipartite structure: its vertices
		may be colored in black and white in such a way that each edge would connect 
		vertices of opposite colors. Namely, black vertices of ${\cal M}$ are the 
		points $x\in f^{-1}(0)$, and white vertices of ${\cal M}$ are the points 
		$x\in f^{-1}(1)$, the vertex degrees being equal to the multiplicities of 
		the corresponding preimages.
\item[\rm 3.]	Inside each face, there is a unique pole of $f$ whose multiplicity
		is equal to the degree of the face. Here the\/ {\em degree of a face} is
		defined as\/ {\em a half} of the number of surrounding edges. We call
		this pole the\/ {\em center} of the face in question.
\end{itemize}
In the opposite direction, if $\cal M$ is a bicolored plane map then:
\begin{itemize}
\item[\rm 4.]	There exists a Belyi function $f$ such that $\cal M$ can be realized
		as a preimage ${\cal M}=f^{-1}([0,1])$.
\item[\rm 5.]	This function $f$ is unique, up to an affine change of the variable $x$.
\item[\rm 6.]	There is a uniquely defined number field $K$ corresponding to $\cal M$ which is
		called the\/ {\em field of moduli} of $\cal M$. The function $f$ can be realized 
		over a number field $L\supseteq K$.
\end{itemize}
\end{theorem}

Statements 4 and 5 represent a particular case of Riemann's existence theorem.
Statement 6 follows from the rigidity of the ramified covering $f:\Cbar\to\Cbar$
and from some general facts of the Galois theory.

The above theorem, being applied to the DZ-pairs, gives the following statement 
(see more details in \cite{PakZvo.I-14}).

\begin{proposition}[DZ-pairs and Belyi functions]\label{prop:DZ-Belyi} 
A pair of complex polynomials $(P,Q)$ is a DZ-pair with a passport $(\al,\be)$ 
if and only if the rational function $f=P/R$, where $R=P-Q$, is a Belyi function 
for a bicolored plane map\/ $\cal M$ with the following characteristics:
\begin{itemize}
\item[\rm 1.]   The map\/ $\cal M$ has $n=\deg P=\deg Q$ edges, $p$ black vertices with
			the degree distribution $\al$, and $q$ white vertices with the degree
			distribution $\be$. {\rm The Euler formula then implies that the number
			of faces is $(n+2)-(p+q)$.}
\item[\rm 2.]   All faces of\/ $\cal M$ except the outer one are of degree~$1$.
\item[\rm 3.]   The number of the faces of\/ $\cal M$ of degree~$1$ is equal to\/ $r=\deg R$.
			In other words, the degree distribution of the faces is equal to $(n-r,1^r)$ where
			$r=(n+1)-(p+q)$.
\end{itemize}		
Furthermore, if $K\subset\Qbar$ is the moduli field of\/ $\cal M$, then it is
possible to find a corresponding DZ-pair such that $P,Q\in K[x]$.
In other words, in this case the realization field $L$ $($see the last statement
of\/ {\rm Theorem~\ref{th:Riemann}}$)$ coincides with the field of moduli $K$.
In particular, an equivalence class of the pair $(P,Q)$ is defined
over\/ $\Q$ if and only if the field of moduli of the map $\cal M$ is $K=\Q$.
\end{proposition}

The characteristic which distinguishes the maps corresponding to DZ-pairs from 
other maps is property 2 of the above theorem. 



\subsection{Weighted trees}
\label{sec:trees}

\begin{figure}[htbp]
\begin{center}
\epsfig{file=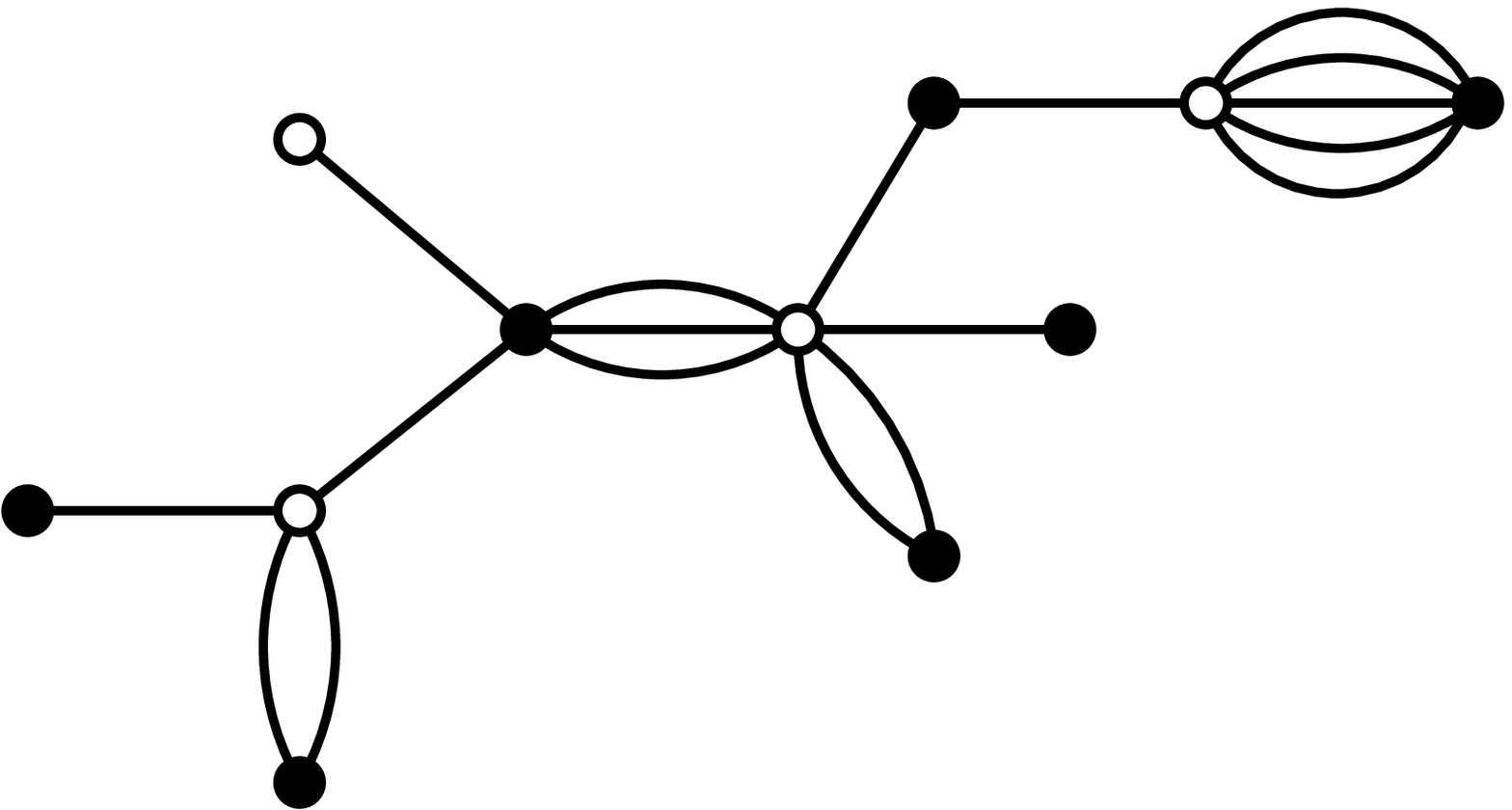,width=5.6cm}
\hspace{0.5cm}
\epsfig{file=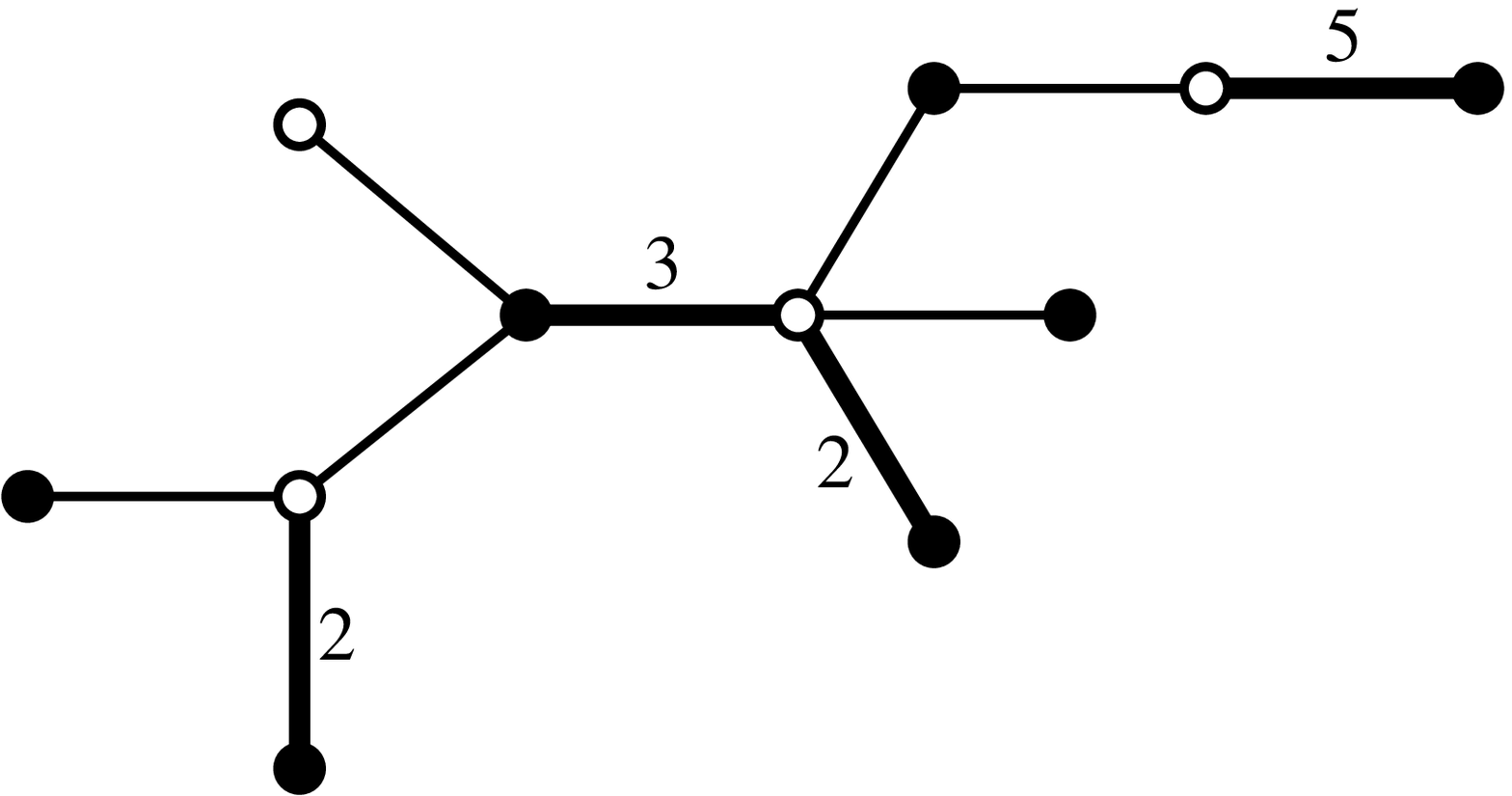,width=5.4cm}
\caption{The passage from a map with all its inner faces being 
of degree~1, to a weighted tree. The weights which are not explicitly 
indicated are equal to~1; the edges of the weight greater than~1 are
drawn thick.}
\label{fig:map->tree}
\end{center}
\end{figure}

We will call the faces other than the outer one {\em inner faces}. The maps
whose all inner faces are of degree 1 can be easily represented in the form 
of weighted trees: just merge every sheaf of parallel edges into one edge and
indicate the number of edges merged together as the {\em weight}\/ of the 
corresponding edge of the weighted tree: see Fig.~\ref{fig:map->tree}. 
Weighted trees are easier to work with than maps.

\begin{definition}[Weighted tree]\label{def:weighted}
A {\em weighted bicolored plane tree}, or a {\em weigh\-ted tree}, or just a
{\em tree}\/ for short, is a bicolored plane tree whose edges are endowed 
with positive integral {\em weights}. The sum of the weights of the edges 
of a tree is called the {\em total weight}\/ or the {\em degree}\/ of the tree. 

The {\em degree}\/ of a vertex is the sum of the weights of the edges incident to this 
vertex. Obviously, the sum of the degrees of black vertices, as well as the sum of the
degrees of white vertices, is equal to the total weight $n$ of the tree. Let the tree
have $p$ black vertices, of degrees $\al_1,\ldots,\al_p$, and $q$ white vertices,
of degrees $\be_1,\ldots,\be_q$, respectively. Then the pair of partitions $(\al,\be)$
of the total weight $n$ of the tree is called its {\em passport}.

Forgetting the weights and considering only the underlying plane tree, 
we speak of a {\em topological tree}. Weighted trees, all of whose edges
are of weight~1, will be called {\em ordinary trees}. Belyi functions for
ordinary trees are polynomials (with the only pole at infinity); they are
usually called {\em Shabat polynomials}.

We call a {\em leaf}\/ a vertex which has only one edge incident to it,
whatever is the weight of this edge. By abuse of language, we will
also call this edge itself a leaf.
\end{definition}

The adjective {\em plane}\/ in the above definition means that the cyclic order 
of branches around each vertex of the tree is fixed, and changing this order will 
in general produce a different plane tree (though the tree considered as a mere 
graph, without ``planar'' structure, remains the same). {\em All trees considered 
in this paper will be endowed with the planar structure}\/; therefore, the adjective 
``plane'' will often be omitted.

The filed of moduli of a unitree is $\Q$, see, e.\,g., \cite{PakZvo.I-14}. Therefore, 
the second part of Theorem~\ref{prop:DZ-Belyi} implies the following statement.

\begin{proposition}[Unitree implies $\Q$]\label{uni->Q}
If a weighted bicolored plane tree is a unitree, then the corresponding equivalence
class of DZ-pairs is defined over\/~$\Q$.
\end{proposition}


\begin{example}[Example \ref{ex:pair-T} revisited]\label{ex:pair-T.bis}
Let us consider the tree shown in Fig.~\ref{fig:T}. It has eight black vertices of degree~3
and twelve white vertices of degree~2, so that its total weight (or degree) is $24$. 
Accordingly, $n=24$, and $\al$ and $\be$ are the following two partitions of $24$: 
$$
\al=(3,3,3,3,3,3,3,3)=3^8, \qquad  \be=(2,2,2,2,2,2,2,2,2,2,2,2)=2^{12}.
$$ 

\begin{figure}[htbp]
\begin{center}
\epsfig{file=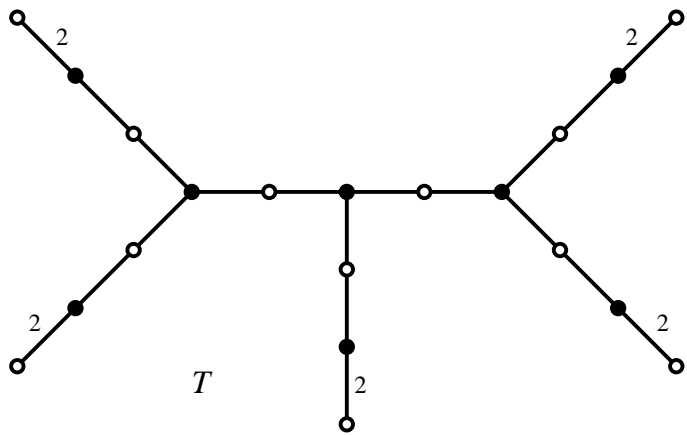,width=9cm}
\end{center}
\caption{One of the sporadic trees of our classification of unitrees
we will speak about further. It is denoted as tree $T$.}\label{fig:T}
\end{figure}

In the corresponding DZ-pair, the polynomial $P$ must have eight roots of
multiplicity~3, the polynomial $Q$ must have twelve roots of multiplicity~2.
In other words, $P=A^3$ with $\deg A=8$, and $Q=B^2$ with $\deg B=12$.
The difference $R=P-Q$ must be of degree $(24+1)-(8+12)=5$.

The general results formulated up to now, being applied to this particular tree, 
imply the following statements:
\begin{itemize}
\item   The mere existence of such a tree implies the existence of polyno\-mials
		  with needed properties.
\item   The fact that there exist polynomials $P$ and $Q$ {\em with rational coefficients}\/ 
		  is a consequence of the fact that there exists a unique tree with the passport $(3^8,2^{12})$.
\end{itemize}
All this can be affirmed without any computations, just by looking at the picture. As to the 
polynomials themselves, they are given in Example~\ref{ex:pair-T}. 
\end{example}


\subsection{Reciprocal polynomials}
\label{sec:reciprocal}

It turns out that technically it is often much more convenient to work not with the
polynomials appearing in DZ-pairs but with their {\em reciprocals}.

\begin{definition}[Reciprocal polynomial]\label{def:reciprocal}
For a polynomial $P$ of degree~$n$, its {\em reciprocal}\/ is 
$P^*(x)=x^n\cdot P(1/x)$.
\end{definition}

In many examples, the reciprocals of polynomials forming a DZ-pair take the 
form of initial segments of power series of some special functions. After having observed 
this phenomenon we learned that it was (re)discovered many times, notably in \cite{Danilov-89}, \cite{Adrianov-97}, \cite{BeuSte-10}.

Assume that  polynomials $P$ and $Q$ form a DZ-pair, so that 
\begin{equation} \label{dere0}
\deg\,(P-Q) \eq (n+1)-(p+q) \eq n-(p+q-1),
\end{equation}
and denote by $m$ the number of edges of the corresponding topological tree.
This tree has $p+q$ vertices, therefore it has $m=p+q-1$ edges.
Considering $P$ and $Q$ as power series we may write condition \eqref{dere0} as 
\begin{equation} \label{dere}
P-Q \eq O(x^{n-m}) \quad \mbox{when} \quad x\to\infty.
\end{equation}
For the reciprocal polynomials condition \eqref{dere0} is transformed into the following one:
\begin{equation} \label{reci0} 
P^* - Q^* \eq x^m\cdot S, 
\end{equation}
where $S$ is a polynomial, or, equivalently, to the condition 
\begin{equation} \label{reci} 
P^* - Q^* \eq O(x^{m}) \quad \mbox{when} \quad x\to 0. 
\end{equation}
For instance, in the Example~\ref{ex:pair-T} the polynomials reciprocal
to \eqref{P-T} and~\eqref{Q-T} and to their difference look as follows:
\begin{eqnarray*}
P^* & = & (1+84x^2+176x^3+2366x^4+13\,536x^5+26\,884x^6 \\
    &   &   +\,\,218\,864x^7+268\,777x^8)^3, \\
Q^* & = & (1+126x^2+264x^3+6195x^4+31\,392x^5+163\,956x^6 \\
    &   &   +\,\,1\,260\,528x^7+3\,531\,639x^8+19\,770\,400x^9 \\
    &   &   +\,\,62\,912\,622x^{10}+94\,024\,776x^{11}+291\,742\,453x^{12})^2, \\
P^*-Q^* & = & x^{19}\,\times\, -2^{38}\cdot 3^3\,(1+62x^2+148x^3+1001x^4+8852x^5).
\end{eqnarray*}


\subsection{Remarks about computation}
\label{sec:computation}

The computation of Belyi functions has recently become a vast domain of research.
A remarkable overview of this activity may be found in \cite{SijVoi-14}, a paper
of 57~pages, with a bibliography of 176~titles. Beside a direct approach, involving
the solution of a system of polynomial equations, the authors of \cite{SijVoi-14}
also discuss complex analytic methods, modular forms methods, and $p$-adic methods.

In order to get an idea of the level of difficulty of such a computation let us return
once again to Example~\ref{ex:pair-T}. A naive approach would be to write 
down polynomials $A \eq \sum_{i=0}^8 u_i x^i$ and $B \eq \sum_{j=0}^{12} v_j x^j$ 
with indeterminate coefficients $u_i$ and $v_j$, and then equate to zero the coefficients 
of degrees from 6 to 24 of the difference $R=A^3-B^2$ . In this way we get 
a system of $24-5=19$ algebraic equations for $9+13=22$ unknowns. Then we may set, for 
example, $u_8=1$, $u_7=0$, and $v_{12}=1$. The system thus obtained (19~equations
with 19~unknowns) will be of degree $25\,509\,168$\,! Obviously, this is not a clever 
way to proceed. 

By the way, the solution we are looking for is unique; all the other solutions of 
this enormous system are ``parasitic'' ones. For example, the system does not
give us any guarantee that the polynomials $A$ and $B$ obtained as its solution 
will be coprime. This condition should be added to the system, but this addition
will make our situation even worse.

Notice, however, that, once the result is obtained, its verification is trivial.

\msk

Taking into account the above considerations, we would like to underline one aspect
of our work: though we do compute Belyi functions for certain individual dessins, 
the most interesting part of the paper is the computation of Belyi functions for
{\em infinite series}\/ of dessins which depend on one or several parameters. 
For infinite series the situation is significantly more complicated than for
individual dessins. Usually, the first thing to do is to compute quite a few 
particular cases, sometimes dozens of them (or to use other heuristics whenever
possible). Then, we need to guess a general pattern of corresponding Belyi functions. 
And, finally, instead of a trivial verification step which was applicable to 
individual dessins, we should provide a {\em proof}, which may turn out to be 
rather laborious. 

In the present paper we obviously do not expose the first step of the above procedure.
What we do is presenting the final results, that is, the general form of Belyi
functions in question, and then we give the proofs whenever they are necessary.
$$
* \qquad * \qquad *
$$

As it was already said, the unitrees comprise ten infinite series, from $A$ to~$J$, 
and ten sporadic trees, from $K$ to $T$. In the subsequent sections
we do not strictly follow the ``alphabetic'' order of trees since we prefer
to underline the structural properties of Belyi functions in question.
Certain Belyi functions are expressed in terms of Jacobi polynomials;
there are others which lead to interesting differential relations; 
we will also encounter compositions, Pad\'e approximants, an application 
to the Hall conjecture, etc.


\section{Stars and binomial series}\label{sec:A}

\begin{figure}[htbp]
\begin{center}
\epsfig{file=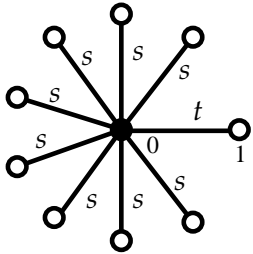,width=3.3cm}
\end{center}
\caption{Star-trees. There are $k$ edges of weight $s$ and one edge
of weight $t$, and ${\rm gcd}(s,t)=1$.}\label{fig:A}
\end{figure}

Our first series, called ``series $A$'' in \cite{PakZvo.I-14}, is composed of 
stars-trees, see Fig.~\ref{fig:A}. All edges except maybe one are of the same weight.
This is a three-parametric series.

Denote the number of leaves of weight $s$ by $k$; then the total weight of the tree
is $n=ks+t$. Clearly, we may put the only black vertex at $x=0$, put
the white vertex of degree $t$ at $x=1$, and assume that both $P$ and $Q$ are monic. 
Then $P(x)=x^n$ and
\begin{eqnarray}\label{Q-A}
Q(x) \eq (x-1)^t\cdot A(x)^s,
\end{eqnarray}
where $A$ is a monic polynomial of degree $k$ whose roots are the white vertices of 
degree~$s$.
Now, condition \eqref{dere} takes the form 
\begin{equation} \label{ert} 
x^n - (x-1)^t\cdot A^s \rav_{x\to \infty}  O(x^{n-(k+1)}).
\end{equation}
The only thing we need to know is the polynomial $A$.

\begin{proposition}
The polynomial $A^*$ reciprocal to $A$ is the initial segment of the binomial 
series for $(1-x)^{-t/s}$ up to the degree $k$:
\begin{eqnarray}\label{A-A}
(1-x)^{-t/s} \rav_{x\to 0}  A^* + O(x^{k+1}).
\end{eqnarray}
\end{proposition}

\paragraph{Proof.} 
Let us pass to reciprocals in \eqref{ert}: we need to obtain $A^*$ such that
$$
1 - (1-x)^t\cdot(A^*)^s  \rav_{x\to 0} O(x^{k+1}).
$$
Let us verify that the polynomial $A^*$ defined in (\ref{A-A}) satisfies
the latter equality. We have:
\begin{equation}\label{xor1} 
A^* \eq (1-x)^{-t/s} + h\cdot x^{k+1},
\end{equation}
where $$h\rav_{x\to 0} O(1).$$ Therefore,
$$
A^*(1-x)^{t/s}  \eq 1+h\cdot x^{k+1}(1-x)^{t/s},
$$
and  
$$
(A^*)^s(1-x)^{t}  \eq \left[1+h\cdot x^{k+1}(1-x)^{t/s}\right]^s \rav_{x\to 0} 
1 + O(x^{k+1})
$$
which concludes the proof.
\hfill $\Box$

\bsk

Some particular cases of formula \eqref{A-A} were previously found
by N.\,Adrianov (unpublished).


\section{Forks and Hall's conjecture}\label{sec:D}

\begin{figure}[htbp]
\begin{center}
\epsfig{file=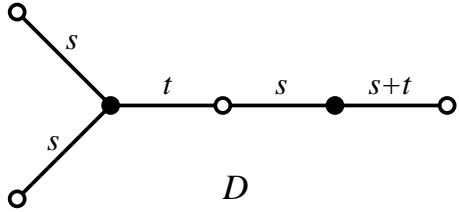,width=6cm}
\end{center}
\caption{Fork-trees. There are exactly two leaves of weight $s$
and exactly one leaf of weigh $s+t$. As usual, ${\rm gcd}(s,t)=1$.}\label{fig:D}
\end{figure}

The two-parametric series of trees shown in Fig.~\ref{fig:D} was called
``series $D$'' in \cite{PakZvo.I-14}.

\subsection{Calculation of DZ-pairs}

This is the only infinite series of unitrees for which we were able to find 
the corresponding DZ-pairs by computer.
Let us introduce the following three {\em quadratic}\/ polynomials:
\begin{center}
\begin{tabular}{lll}
$A$ & -- & the roots of $A$ are two black vertices of degree $2s+t$; \\
$B$ & -- & the roots of $B$ are two white vertices of degree $s+t$; \\
$C$ & -- & the roots of $C$ are two white vertices of degree $s$.
\end{tabular}
\end{center}

\begin{proposition}
We have $P=A^{2s+t}$ and $Q=B^{s+t}\cdot C^s$, where
\begin{eqnarray}
A & = & x^2-(3s+t)(3s+2t); \label{A-D} \\
B & = & x^2-6s\cdot x+(3s-2t)(3s+t); \label{B-D} \\
C & = & x^2+6(s+t)\cdot x+(3s+2t)(3s+5t). \label{C-D}
\end{eqnarray}
\end{proposition}

\paragraph{Proof.} By \eqref{reci}, we must prove that 
\begin{eqnarray}\label{D-Maple}
(A^*)^{2s+t} - (B^*)^{s+t}\cdot(C^*)^s \eq O(x^5).
\end{eqnarray}
Clearly, we may assume that the sum of the roots of $A$ equals zero. 
Write 
$$
A^* = 1-ax^2, \qquad
B^* = 1-b x+cx^2, \qquad
C^* = 1+d x+ex^2,
$$
and calculate, with the help of Maple, the first five coefficients of the 
Taylor series in the left-hand side of (\ref{D-Maple}). Equate now the 
expressions thus obtained to zero and solve the corresponding system in 
the unknowns $a,b,c,d,e$. Maple returns two solutions: 
$$
a=-e, \quad b=0, \quad c=e, \quad d=0, \quad e=e,
$$ 
and 
$$ 
a \eq {\frac {{b}^{2} \left( 9\,{s}^{2}+9\,ts+2\,{t}^{2} \right) }{36{s}^{2}}}, \qquad  
b=b, \qquad 
c={\frac {{b}^{2} \left( 9\,{s}^{2}-3\,ts-2\,{t}^{2} \right) }{36{s}^{2}}},
$$ 
$$
d \eq {\frac { \left( t+s \right) \\ \mbox{}b}{s}}, \qquad
e \eq {\frac {{b}^{2} \left( 9\,{s}^{2}+10\,{t}^{2}+21\,ts \right) }{36{s}^{2}}}.
$$
Rejecting the first solution, for which the roots of $A$, $B$ and $C$ coincide, 
and making an additional normalization by setting the $b=6s$, we obtain 
formulas \eqref{A-D}, \eqref{B-D} and \eqref{C-D}.
\hfill$\Box$


\subsection{An application: Danilov's theorem}

In 1971, M. Hall, Jr. \cite{Hall-71} suggested the following two conjectures. 
\begin{enumerate}
\item	There exists a constant $c$ such that for all positive integers $a,b$, 
		$a^3\neq b^2$, the inequality 
		$$ 
		|a^3-b^2| \,>\, c\cdot a^{1/2} 
		$$ 
		holds. 
\item	The exponent $1/2$ in the above inequality cannot be improved.
		Namely, for every $\eps>0$ there exists a constant $C(\eps)$ such that there are 
		infinitely many pairs of integers $(a,b)$ satisfying the inequality
		$$
		|a^3-b^2| \,\le\, C(\eps)\cdot a^{1/2+\eps}. 
		$$
\end{enumerate}
This first conjecture is neither proved nor disproved. However, a general belief is 
that in order to be true it should be modified as follows: for each $\eps>0$ there 
exists a constant $c(\eps)$ such that for all positive integers $a,b$, $a^3\neq b^2$, 
the inequality 
$$
|a^3-b^2| \,>\, c(\eps)\cdot a^{1/2-\eps} 
$$
holds. In this form the conjecture is a corollary of the famous $ABC$-conjecture 
(see, e.\,g., \cite{Lang-90}, \cite{BeuSte-10} for further details).

As to the second conjecture, in 1982 Danilov \cite{Danilov-82} proved its stronger
version. His result is interesting for us since in his proof he used, in a slightly 
different normalization, the above polynomials $A,B,C$, see \eqref{A-D}, \eqref{B-D},
\eqref{C-D}, with the parameters $s=t=1$. 

\begin{proposition}[Danilov's theorem]\label{prop:Hall}
There exists a constant $C$ such that there are infinitely many
pairs of integers $(a,b)$ satisfying the inequality
\begin{equation} \label{ner}
|a^3-b^2| \,\le\, C\cdot a^{1/2}.
\end{equation}
\end{proposition}

\paragraph{Proof.} Specializing (\ref{A-D}), (\ref{B-D}) and (\ref{C-D}) for
$s=t=1$ and computing the difference $P-Q$ we get
$$
(x^2-20)^3 - (x^2-6x+4)^2(x^2+12x+40) \eq 1728x - 8640.
$$
Substituting $x=2z$ and dividing both parts by 8  we get
\begin{equation} \label{giv}
(2z^2-10)^3 - (2z^2-6z+2)^2(2z^2-12z+20) \eq 432z - 1080.
\end{equation}
Let us now consider the factor $2z^2-12z+20=2(z-3)^2+2$ and try to make it a perfect square;
then \eqref{giv} will give us a relatively ``small'' difference between a cube and a square. 
To do that we have to solve the Diophantine equation
\begin{equation}\label{urrr}
u^2 - 2v^2 = 2,
\end{equation}
where $v=z-3$. 

The last equation is a Pell-like equation, that is an equation of the form 
$$
u^2 - Dv^2 = m,
$$ 
where $D>0$ is a square-free integer and $m\in \Z$. For $m=1$ this equation 
is a usual Pell equation, and it is well known that any  Pell equation has 
infinitely many integer solutions. Pell-like equations not necessarily have 
integer solutions. However, if at least one such solution $(u_0,v_0)$ exists, 
then we can obtain infinitely many solutions $(u_n,v_n)$ using the following
recursion:
$$
u_n+v_n\sqrt{D} \eq (u_{n-1}+v_{n-1}\sqrt{D})(k+l\sqrt{D})
$$
where $(k,l)$ is the minimum solution of the equation $k^2-Dl^2=1$. In our case,
$(k,l)=(3,2)$.

Equation \eqref{urrr} does have an integer solution $(u_0,v_0)=(2,1)$. Returning 
to~\eqref{giv}, it is easy to verify that for all $z\ge 3$ one has
$$
432z-1080 \,<\, 216\sqrt{2}\cdot(2z^2-10)^{1/2},
$$
which proves the theorem: there are infinitely many pairs of integers $(a,b)$ 
satisfying \eqref{ner}, with the constant $C=216\sqrt{2}$. 
\hfill$\Box$

\bsk

The same polynomials $A,B$ with the parameters $s=t=1$ were used by 
Dujella~\cite{Dujella-10} for constructing an infinite series of pairs 
of polynomials $P,Q$ with the following properties: (a) $\deg P=2k$, $\deg Q=3k$; 
(b) $P$ and $Q$ are {\em not}\/ coprime; (c) $\deg(P^3-Q^2)=k+5$, so 
that the minimum degree $k+1$ is not attained, though the discrepancy 
remains bounded; (d) in return, $P$ and $Q$ are defined over $\Q$. 

\ssk 

Using other DZ-pairs, Danilov \cite{Danilov-89} and Beukers and Stewart \cite{BeuSte-10}
obtained results similar to Proposition \ref{prop:Hall} for the differences 
between integer powers $a^n$~and~$b^m$.



\section{Jacobi polynomials}

\subsection{Trees of this section}

\begin{figure}[htbp]
\begin{center}
\epsfig{file=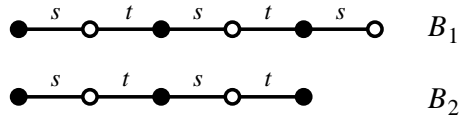,width=6cm}
\end{center}
\caption{Series $B_1$ and $B_2$: chain-trees}\label{fig:B}
\end{figure}

Davenport--Zannier pairs for the series of trees considered in this section are 
expressed in terms of Jacobi polynomials. The trees in question are constructed 
as follows. First, we take chain-trees with alternating edge weights $s,t,s,t,\ldots$, 
see Fig.~\ref{fig:B}. We must distinguish chains of odd and even length since in one 
case both ends are of the same color while in the other case they are of different 
colors.

\begin{figure}[htbp]
\begin{center}
\epsfig{file=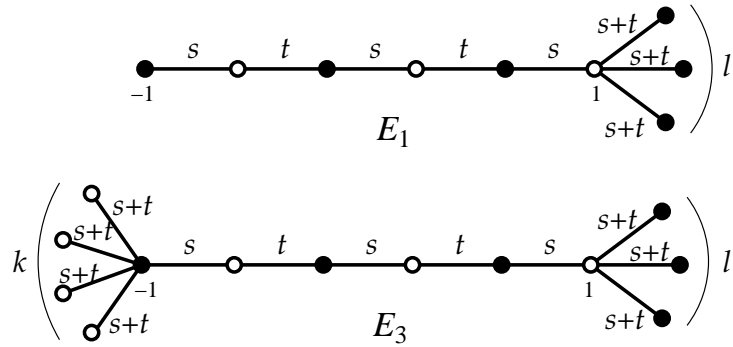,width=9.6cm}
\end{center}
\caption{Series $E_1$ and $E_3$: odd double brushes}\label{fig:E1-E3}
\end{figure}

\begin{figure}[htbp]
\begin{center}
\epsfig{file=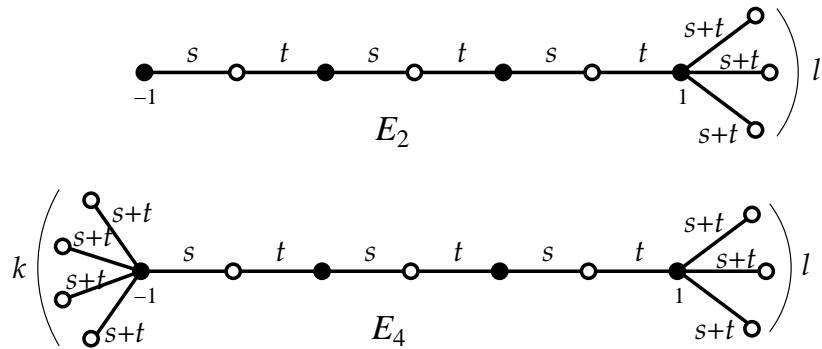,width=10.8cm}
\end{center}
\caption{Series $E_2$ and $E_4$: even double brushes}\label{fig:E2-E4}
\end{figure}

Then, we have a right to attach to the end-points an arbitrary number of leaves
of the weight $s+t$. In this way we obtain ``odd'' series $E_1,E_3$ and ``even''
series $E_2,E_4$, see Figs.~\ref{fig:E1-E3} and \ref{fig:E2-E4}. We call these
series ``double brushes''. Note that any of the parameters $k,l$, and also both
of them, may be equal to zero. Thus, $B_1$ and $E_1$ are particular cases of
$E_3$, and $B_2$ and $E_2$ are particular cases of $E_4$.

\begin{figure}[htbp]
\begin{center}
\epsfig{file=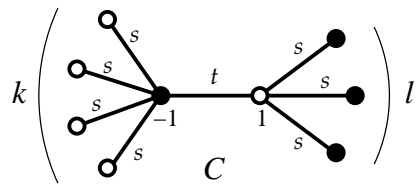,width=5.4cm}
\end{center}
\caption{Series $C$: trees of diameter 3}\label{fig:C}
\end{figure}

There are two exceptions from the above construction. The first is when the chain 
part consists of a single edge, so that there is no alternance of weights.
We thus obtain the series $C$, see Fig.~\ref{fig:C}. In contrast to
the general case, now the weight of leaves may be smaller than the 
weight of the edge between the leaves.

The second exception is when the chain part consists of two edges. In this case 
it is possible to attach exactly one leaf of weight $s+t$ to one of the ends and 
exactly two leaves of weight $s$ (or $t$, to ensure the weight alternance) to the 
other end. In this way, we get the series of forks $D$ already studied in Sect.~\ref{sec:D}.


\subsection{Jacobi polynomials: preliminaries}

Let us recall some general facts concerning Jacobi polynomials; for more advances
and detailed treatment see, for example, \cite{Szego-39} or \cite{AbrSte-72}.

The classical Jacobi polynomials $J_n(a,b,x)$, $\deg J_n=n$, are defined for the
parameters $a,b\in\R$, $a,b>-1$, as orthogonal polynomials with respect to the 
measure on the segment $[-1,1]$, given by the density $(1-x)^a(1+x)^b$. The 
restriction $a,b>-1$ is necessary in order to ensure the integrability.
The polynomial $J_n(a,b,x)$ can also be defined as a unique polynomial solution of 
the differential equation 
\begin{equation} \label{urav}
(1-x^2)y''+\left[b-a-(a+b+2)x\right]y'+n(n+a+b+1)y \eq 0,
\end{equation} 
satisfying the condition $J_n(a,b,1)=\binom{n+a}{n},$ or by the explicit formula 
\begin{equation} \label{form}
J_n(a,b,x) \eq \sum_{k=0}^n\binom{n+a+b+k}{k}\binom{n+a}{n-k}\left(\frac{x-1}{2}\right)^k.
\end{equation} 

Notice that equation \eqref{urav} can be written in the form 
\begin{equation}\label{u}
(1-x^2)Y''+\left[a-b+(a+b-2)x\right]Y'+(n+1)(n+a+b)Y \eq 0,
\end{equation}
where $Y=(1-x)^a(1+x)^b\cdot y$, implying that the function 
\begin{equation}\label{fun}
(1-x)^a(1+x)^b\cdot J_n(a,b,x)
\end{equation} 
satisfies \eqref{u}. 

It follows from 
\eqref{form} that $J_n(a,b,x)$ are also 
polynomials in parameters $a$ and~$b$. Therefore, their definition can be 
extended to arbitrary (even complex) values of these parameters. 
These {\em generalized}\/ Jacobi polynomials still satisfy~\eqref{urav}, 
although they are no longer orthogonal with respect to a measure on the
segment $[-1,1]$. 
Similarly, since the function \eqref{fun} may be represented as a power series in $x$ 
whose coefficients are polynomials in $a,b$, this function satisfies equation \eqref{u} 
for arbitrary $a$ and $b$.

The following key observation will be used in subsequent proofs. If, in the
differential operator \eqref{urav}, we replace $n$ with $n+a+b$, $a$ with $-a$,
and $b$ with $-b$, we get exactly the differential operator \eqref{u}.
Therefore, $J_{n+a+b}(-a,-b,x)$ along with \eqref{fun} satisfies \eqref{u}. 
The last statement, however, should be taken with caution: the subscript $n+a+b$ 
must be a non-negative integer since it is the degree of a polynomial.

Notice that if $a$ and $b$ do not satisfy the inequalities $a,b>-1$, then 
the degree in $x$ of the polynomial $J_n(a,b,x)$ defined by \eqref{form} may drop 
down below~$n$. Indeed, \eqref{form} implies that the leading coefficient of 
$J_n(a,b,x)$ is equal to
\begin{eqnarray}\label{eq:jacobi}
\frac{1}{2^n}\binom{2n+a+b}{n} \eq \frac{1}{2^n\cdot n!}\prod_{i=n+1}^{2n} (a+b+i)\,.
\end{eqnarray}
Hence, in order to obtain a polynomial of degree $n$ we must require that the sum 
$a+b$\/ does not take values $-(n+1)$, $-(n+2)$, \ldots, $-2n$. In particular, 
this is always true if $a$ and $b$ are real and $n\ge -(a+b)$ or, equivalently, 
$n+a+b\ge 0$.

Along with the density $(1-x)^a(1+x)^b$, which is defined on $[-1,1]$, we will use 
the multivalued complex function $(z-1)^a(z+1)^b$ (note the change of the sign of 
the term in the first parenthesis). Clearly, this function has three ramification 
points $-1,1,\infty$. Further, observe that 
if $a+b\in \Z$, then any germ of $(z-1)^a(z+1)^b$ defined near a non-singular 
point $z_0$ extends to a function $\mu(z)$ which is single-valued in any domain $U$ 
obtained from $\C\P^1$ by removing a simple curve connecting $-1$ and $1$. 
Indeed, in such $U$ the function $\mu(z)$ may have a ramification only at infinity. 
On the other hand, since the analytic continuation of $\mu(z)$ along a loop around 
infinity is $e^{2\pi (a+b)i}\mu(z)$, we see that $\infty$ is not a ramification 
point since $a+b\in \Z$. In particular,  $\mu(z)$ can be expanded into a Laurent 
series at infinity,
$$
\mu(z) \eq c_{a+b}z^{a+b}+c_{a+b-1}z^{a+b-1} + \ldots\ \  .
$$ 
Finally, if $a$ and $b$ are rational numbers, say
\begin{equation} \label{para}
a \eq \frac{n_1}{m}, \quad b\eq \frac{n_2}{m}, \quad n_1,n_2,m\in \Z,
\end{equation} 
then any $\mu(z)$ as above satisfies the condition 
$$
\mu(z)^m \eq (z-1)^{n_1}(z+1)^{n_2},
$$ 
implying that $\mu(z)$ is defined up to a multiplication by an $m$th root of unity, 
and that for a certain choice of this root  the equality $c_{a+b}=1$ holds. By abuse
of notation, below we will always use the expression $(z-1)^a(z+1)^b$ to denote 
the function $\mu(z)$ which satisfies the equality  $c_{a+b}=1$.

\begin{lemma}\label{lemma} 
Assume that $a$ and $b$ are rational numbers which satisfy the condition 
$a+b\in \Z$. Then for any $n\geq -(a+b)$ the equality 
\begin{equation}\label{jacobi} 
\left(\frac{z-1}{2}\right)^a\left(\frac{z+1}{2}\right)^b
J_n(a,b,z)-J_{n+a+b}(-a,-b,z)\, \rav_{z\to \infty}\, O(z^{-(n+1)})
\end{equation} 
holds. 
\end{lemma}

\paragraph{Proof.} As it was mentioned above, the function \eqref{fun} satisfies 
the differential equation \eqref{u}, where the function $\nu(x)=(1-x)^a(1+x)^b$ 
is assumed to be defined on $[-1,1]$. However, since this function is analytic 
near the origin, we can consider its analytic continuation $\nu(z)$, and the 
function $\nu(z)J_n(a,b,z)$ will satisfy \eqref{u} in the domain $U$ as above. 
Furthermore, if \eqref{para} holds, then 
$$
\nu(z)^m \eq (-1)^{n_1}\Big((z-1)^a(z+1)^b\Big)^m,
$$ 
implying that the function $(z-1)^a(z+1)^bJ_n(a,b,z)$ also satisfies \eqref{u} 
in $U$.

Since  the polynomial $J_n(a,b,x)$ satisfies the differential equation \eqref{urav}, 
we conclude that the functions 
$$
Y_1 \eq \left(\frac{z-1}{2}\right)^a\left(\frac{z+1}{2}\right)^bJ_n(a,b,z) \quad  
{\rm and} \quad  Y_2 \eq J_{n+a+b}(-a,-b,z)
$$ 
both satisfy the differential equation 
\begin{equation}\label{difeq} 
L_{n}^{a,b}(Y) \eq 0,
\end{equation} 
where 
$$
L_{n}^{a,b} \eq (1-z^2)\frac{d^2}{dz^2}+\left[a-b+(a+b-2)z\right]\frac{d}{d z}+(n+1)(n+a+b).
$$ 
This implies that the function $Y_0=Y_1-Y_2$ also satisfies this equation. 
On the other hand, it is easy to see that if $Y(z)$ is a function whose Laurent expansion 
at infinity is 
$$
Y \eq C_dz^d+C_{d-1}z^{d-1}+\ldots ,
$$ 
then 
$$
L_n^{a,b}(Y) \eq \widetilde C_dz^d+\widetilde C_{d-1}z^{d-1}+\ldots
$$
where 
\begin{eqnarray*}
\widetilde C_d & = & -d(d-1)+d(a+b-2)+(n+1)(n+a+b) \\
               & = & (n+a+b-d)(d+n+1).
\end{eqnarray*}
Therefore, if $Y$ satisfies \eqref{difeq} and $C_d\neq 0$ while $\widetilde C_d = 0$, 
we should have either $d=n+a+b$ or $d=-(n+1)$. Finally, \eqref{form} implies that 
the leading terms of both $Y_1$ and $Y_2$ are equal to
$$
\frac{1}{2^{n+a+b}}\binom{2n+a+b}{n}z^{n+a+b}. 
$$ 
Therefore, the degree of the leading term of their difference $Y_0=Y_1-Y_2$ is less than
$n+a+b$, hence the only possible case is $d=-(n+1)$, implying \eqref{jacobi}.
\hfill$\Box$


\subsection{Double brushes of even length}\label{sec:E-even}

Let $\cal T$ be a weighted tree from the series $E_4$ or of its two particular 
cases $E_2$ or~$B_2$, see Figs. \ref{fig:E2-E4} and \ref{fig:B}. Denote by $r$ 
the number of white vertices of $\cal T$ which are not leaves. Then the total 
weight of $\cal T$ is equal to $(s+t)(k+l+r)$ and  the total number of edges 
is equal to $k+l+2r$. Clearly,
\begin{eqnarray}
P & = & (x-1)^{l(s+t)+t}(x+1)^{k(s+t)+s}\cdot A^{s+t}, 	\label{E4-PA} \\ 
Q & = & B^{s+t}								\label{E4-PB}
\end{eqnarray}
for some polynomials $A$ and $B$ with $\deg A=r-1$, $\deg B=k+l+r$. 
Furthermore, by \eqref{dere}, we must have:
$$
P-Q\, \rav_{x\to \infty}\, O(x^{m}),
$$ 
where 
\begin{equation}\label{fvalue} 
m \eq (s+t)(k+l+r)-(k+l+2r) \eq (k+l+r)(s+t-1)-r.
\end{equation}

\begin{proposition} \label{pro} 
The polynomials $P$ and $Q$ may be represented as follows:
\begin{eqnarray}\label{P-E2-E4}
P(x) \,=\, \left(\frac{x-1}{2}\right)^{l(s+t)+t}\cdot \left(\frac{x+1}{2}\right)^{k(s+t)+s}\cdot
J_{r-1}(a,b,x)^{s+t},
\end{eqnarray}
where $J_{r-1}(a,b,x)$ is the Jacobi polynomial with parameters
\begin{eqnarray}\label{param-E2-E4}
a \eq \frac{l(s+t)+t}{s+t} \qquad \mbox{and} \qquad b \eq \frac{k(s+t)+s}{s+t},
\end{eqnarray}
and
\begin{eqnarray}\label{Q-E2-E4}
Q(x) \,=\, J_{k+l+r}(-a,-b,x)^{s+t}.
\end{eqnarray}
\end{proposition}

\paragraph{Proof.} 
Since the polynomials $A$ and $B$ in \eqref{E4-PA}, \eqref{E4-PB} are defined 
in a unique way up to a multiplication by a scalar factor, it is enough to show that 

{\small
\begin{equation}\label{formula} 
\left(\frac{x-1}{2}\right)^{l(s+t)+t}
\left(\frac{x+1}{2}\right)^{k(s+t)+s}
J_{r-1}(a,b,x)^{s+t} - J_{k+l+r}(-a,-b,x)^{s+t}\, \rav_{x\to \infty}\, O(x^{m}),
\end{equation}
}
where $a$ and $b$ are given by \eqref{param-E2-E4}, and $m$, by \eqref{fvalue}.

Represent the left side of \eqref{formula} as a product of two factors using the formula 
\begin{equation}\label{produ}
u^{s+t}-v^{s+t} \eq (u-v)(u^{s+t-1}+u^{s+t-2}v+\dots + v^{s+t-1}),
\end{equation}
where 
$$
u \eq \left(\frac{x-1}{2}\right)^a\left(\frac{x+1}{2}\right)^bJ_{r-1}(a,b,x), \qquad 
v \eq J_{k+l+r}(-a,-b,x),
$$
It is easy to see that both $u$ and $v$ are $O(x^{k+l+r})$ near infinity.
Let us consider the difference $u-v.$
Clearly, 
$$
k+l+r \eq r-1+a+b.
$$ 
Furthermore, since $k,l,r\geq 0$ the inequality 
$$
r-1\, \geq\, -(a+b) \eq -(k+l+1)
$$ 
holds. Therefore, by Lemma \ref{lemma}, we have: 
$$
u-v\,\rav_{x\to \infty}\,O(x^{-r}).
$$ 
On the other hand, 
$$
u^{s+t-1}+u^{s+t-2}v+\dots + v^{s+t-1}\, \rav_{x\to \infty}\, O(x^{(k+l+r)(s+t-1)}).
$$
Thus, 
$$
u^{s+t}-v^{s+t}\,\rav_{x\to \infty}\, O(x^{m})
$$
as required.
\hfill$\Box$

\begin{remark}
Belyi functions for the series $E_2$ and $E_4$ with the parameters $s=t=1$ 
were first calculated in the thesis of Nicolas Magot in 1997 \cite{Magot-97}. 
A different proof, proposed by Don Zagier, was given in Ch.~2 of \cite{LanZvo-04}. 
We used Zagier's proof as a model for the above construction.
\end{remark}


\subsection{Series $E_1$ and $E_3$: double brushes of odd length}\label{sec:E-odd}

Let now $\cal T$ be  a weighted tree from the series $E_3$ or of its two
particular cases $E_1$ and $B_1$, see Figs. \ref{fig:E1-E3} and \ref{fig:B}. 
As above, denote by $r$ the number of white vertices of $\cal T$ 
which are not leaves, so that the total weight of $\cal T$ is $(s+t)(k+l+r)+s$,  
and the total number of edges is $k+l+2r+1$. Now we must find polynomials 
$P$ and $Q$ such that 
\begin{eqnarray}
P & = & (x+1)^{k(s+t)+s}\cdot A^{s+t}, \label{E3-PA} \\
Q & = & (x-1)^{l(s+t)+s}\cdot B^{s+t} \label{E3-QB}
\end{eqnarray}
for some polynomials $A$ and $B$ with $\deg A = l+r$ and $\deg B = k+r$, and  
$$
P-Q\, \rav_{x\to \infty}\, O(x^{m}),
$$ 
where 
\begin{equation}\label{fvalue1} 
m \eq (s+t)(k+l+r)+s-(k+l+2r+1) \eq (k+l+r)(s+t-1)+s-r-1.
\end{equation}

\begin{proposition}
The polynomials $P$ and $Q$ may be represented as follows:
\begin{eqnarray}\label{P-E1-E3}
P(x) \,=\,  \left(\frac{x+1}{2}\right)^{k(s+t)+s}\cdot J_{l+r}(a,b,x)^{s+t},
\end{eqnarray}
where $J_{l+r}(a,b,x)$ is the Jacobi polynomial with the parameters
\begin{eqnarray}\label{param-E1-E3}
a \eq -\frac{l(s+t)+s}{s+t} \qquad \mbox{and} \qquad b \eq \frac{k(s+t)+s}{s+t},
\end{eqnarray}
and
\begin{eqnarray}\label{Q-E1-E3}
Q(x) \,=\, \left(\frac{x-1}{2}\right)^{l(s+t)+s}\cdot J_{k+r}(-a,-b,x)^{s+t}.
\end{eqnarray}
\end{proposition}

\paragraph{Proof.} We must show that 
{\small
\begin{equation}\label{formula2} 
\left(\frac{x+1}{2}\right)^{k(s+t)+s}J_{l+r}(a,b,x)^{s+t} -  
\left(\frac{x-1}{2}\right)^{l(s+t)+s}J_{k+r}(-a,-b,x)^{s+t}\, \rav_{x\to \infty}\, O(x^{m})
\end{equation}
}
where 
$$
a \eq -\frac{l(s+t)+s}{s+t}, \qquad b \eq \frac{k(s+t)+s}{s+t},
$$
and $m$ is defined by \eqref{fvalue1}.

Equality \eqref{formula2} is equivalent to the equality 
{\small
\begin{equation}\label{formula3} 
\left(\frac{x-1}{2}\right)^{-(l(s+t)+s)}\left(\frac{x+1}{2}\right)^{k(s+t)+s}
J_{l+r}(a,b,x)^{s+t} - J_{k+r}(-a,-b,x)^{s+t} \eq O(x^{p}),
\end{equation}
}
where 
$$
p \eq m-(l(s+t)+s) \eq (k+r)(s+t-1) - (l+r+1).
$$
On the other hand, since 
$$
k+r \eq (l+r)+a+b
$$ 
and 
$$
l+r\, \geq\, -(a+b) \eq l-k,
$$
it follows from Lemma \ref{lemma} that 
$$
\left(\frac{x-1}{2}\right)^a\left(\frac{x+1}{2}\right)^bJ_{l+r}(a,b,x)-J_{k+r}(-a,-b,x) 
\eq O(x^{-(l+r+1}),
$$
implying in the same way as in Proposition \ref{pro} that \eqref{formula3} holds.
\hfill$\Box$


\subsection{Series $C$ and $B$}\label{sec:C}

The series $C$ is a particular case of the series $E$ of odd length corresponding 
to the case of $r$ equal to zero. In order to adjust the notation (which is slightly 
different for the series $E$ and $C$) we must set $r=0$ and change $s$ to $t$ and 
$t$ to $s-t$ in formulas \eqref{P-E1-E3}--\eqref{Q-E1-E3}. Thus, 
\begin{eqnarray}\label{P-C}
P(x) \,=\,  \left(\frac{x+1}{2}\right)^{ks+t}\cdot J_l(a,b,x)^{s},
\end{eqnarray}
where $J_l(a,b,x)$ is the Jacobi polynomial of degree $l$ with parameters 
\begin{eqnarray}\label{param-C}
a \eq -\frac{ls+t}{s} \qquad \mbox{and} \qquad b \eq \frac{ks+t}{s},
\end{eqnarray}
while
\begin{eqnarray}\label{Q-C}
Q(x) \,=\, \left(\frac{x-1}{2}\right)^{ls+t}\cdot J_k(-a,-b,x)^{s}.
\end{eqnarray}
Finally, it is clear that the series $B_1$ and $B_2$ (chains of odd and even length) 
are particular cases of the series $E_3$ and $E_4$, so that the Davenport--Zannier 
pairs for $B_1$ and $B_2$ are obtained from those for $E_3$ and $E_4$ by setting $k=l=0.$


\subsection{Pad\'e approximants}\label{sec:pade}

The above results can be interpreted in terms of Pad\'e approximants for 
the function $(1-x)^a(1+x)^b$. Recall that if 
$$
f(x) \eq \sum_{k=0}^{\infty}c_kx^k
$$ 
is a formal power series, then its {\em Pad\'e approximant}\/ of order $[n/m]$ 
at zero is a rational function $p_n(x)/q_m(x)$, where $p_n(x)$ is a polynomial of 
degree $\leq n$ and  $q_m(x)$ is a polynomial of degree $\leq m$, such that 
\begin{equation}\label{p1}
f(x)-\frac{p_n(x)}{q_m(x)}\, \rav_{x\to 0}\, O(x^{n+m+1}).
\end{equation}
Defined in this way, Pad\'e approximants do not necessarily exist. However, if an 
approximant of a given order exists, it is unique.

Linearizing the problem by requiring that 
\begin{equation}\label{p2}
q_m(x)f(x) - p_n(x)\, \rav_{x\to 0}\, O(x^{n+m+1})
\end{equation}
we arrive to the notion of a {\em Pad\'e form}\/ $(p_n,q_m)$ of order $[n/m]$.
Being defined by linear equations, Pad\'e forms always exist (in general, 
\eqref{p2} does not imply \eqref{p1} since $q_m(x)$ may vanish at zero), 
and the Pad\'e form of a given order is defined in a unique way up 
to a multiplication by a constant.

Keeping the notation of Sect.~\ref{sec:E-even} we may now reformulate the 
condition for $P$ and $Q$ to be a Davenport--Zannier pair for the series $E$ 
of even length as follows (a similar result is also true for the series $E$ of 
odd length).
\begin{proposition}[Pad\'e forms, even case]\label{prop:pade-even}
Let polynomials $A$ and $B$ be like in formulas \eqref{E4-PA} and \eqref{E4-PB}.
Then the pair of their reciprocals $(A^*,B^*)$ is the Pad\'e form of order 
$[r-1/k+l+r]$ for the function  $(1-x)^a(1+x)^b$ with parameters
\begin{eqnarray}
a \eq \frac{l(s+t)+t}{s+t} \qquad \mbox{and} \qquad b \eq \frac{k(s+t)+s}{s+t}.
\end{eqnarray}
\end{proposition}

\paragraph{Proof.} Since the pairs $(P,Q)$ and $(A,B)$ are both defined up 
to a multiplication by a constant, it is enough to show that 
\begin{equation}\label{nado1}
(x-1)^{l(s+t)+t}(x+1)^{k(s+t)+s}\cdot A^{s+t}-B^{s+t}\,	\rav_{x\to \infty}\, O(x^{p}),
\end{equation}
where 
$$
p \eq (k+l+r)(s+t-1)-r.
$$
By definition of Pad\'e forms we have: 
$$
(1-x)^a(1+x)^bA^*- B^*\, \rav_{x\to 0}\, O(x^{k+l+2r}),
$$ 
implying that 
\begin{equation}\label{nado2}
(1-x)^{l(s+t)+t}(1+x)^{k(s+t)+s}\cdot (A^*)^{s+t}-(B^*)^{s+t}\,	\rav_{x\to 0}\, O(x^{k+l+2r}),
\end{equation}
(here we use formula \eqref{produ} again though now the factors involved are series 
by non-negative powers of $x$). Finally, substituting $1/x$ in place of $x$ in \eqref{nado2} 
and multiplying both sides by 
$$
x^{(k+l+r)(s+t)}
$$ 
we obtain \eqref{nado1}.
\hfill$\Box$

\begin{proposition}[Pad\'e forms, odd case]\label{prop:pade-odd}
Let polynomials $A$ and $B$ be like in formulas \eqref{E3-PA} and \eqref{E3-QB}.
Then the pair of their reciprocals $(A^*,B^*)$ is the Pad\'e form of order 
$[l+r/k+r]$ for the function  $(1-x)^a(1+x)^b$ with parameters
\begin{eqnarray}
a \eq -\frac{l(s+t)+t}{s+t} \qquad \mbox{and} \qquad b \eq \frac{k(s+t)+s}{s+t}.
\end{eqnarray}
\end{proposition}

The proof is similar to the previous one, so we omit it.

\begin{remark}[On Pad\'e approximants]\label{rem:pade}
From\, the\, computational\, point\, of\, view, a great advantage of Pad\'e approximants 
is due to the fact that the equations describing them are linear. This observation
remains true even in the case like ours when the polynomials in question are
known explicitly. One has to use some astute tricks in order to make Maple work 
with Jacobi polynomials whose parameters do not satisfy the condition $a,b>-1$. 
At the same time, the computation of Pad\'e approximants is instantaneous.


A vast literature is devoted to the study of Pad\'e approximants for some
particular functions. This is the case, for example, for the exponential
function. To our surprise, we did not find any research concerning Pad\'e
approximants for the function $(1-x)^a(1+x)^b$. By the way, our Lemma~\ref{lemma}
can also be reformulated as a result about Pad\'e forms for this function.
\end{remark}


\section{Series $F$ and $G$: trees of diameter 4}

\begin{figure}[htbp]
\begin{center}
\epsfig{file=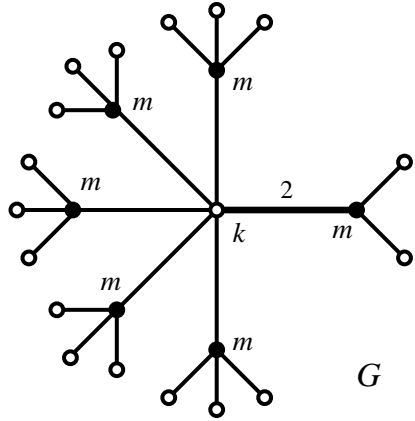,width=5.4cm}
\end{center}
\caption{Series $G$. The degree of the central vertex is $k$, the number of 
branches (and the number of black vertices) is $k-1$.}\label{fig:G}
\end{figure}

\begin{figure}[htbp]
\begin{center}
\epsfig{file=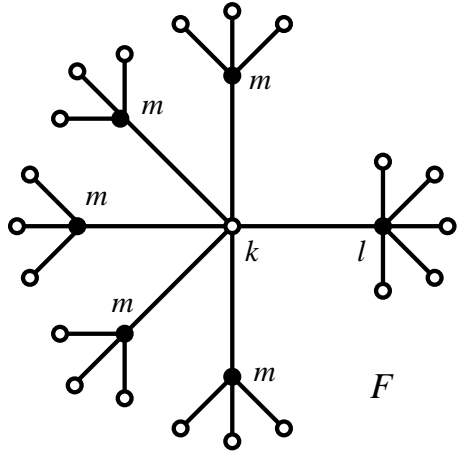,width=6cm}
\end{center}
\caption{Series $F$}\label{fig:F}
\end{figure}

Below we find DZ-pairs for the series $F$ and $G$, see Figs.~\ref{fig:F} and \ref{fig:G}, 
using their relations to differential equations. For the series $F$, which consists of 
ordinary trees, the corresponding formulas are particular cases of the formulas for 
{\em Shabat polynomials}\/ for trees of diameter four, first calculated by 
Adrianov~\cite{Adrianov}. 

Since any tree from the series $F$ is ordinary, the degree of $R=P-Q$ is zero, 
that is $R=c$ for some $c\in \C.$ Therefore, in order to describe the corresponding 
DZ-pair it is enough to find $P$ and $c$. This is equivalent to the finding of the  
Shabat polynomial corresponding to the tree. 
Similarly, for trees from the series $G$ the degree of $R$ is one, and 
it is technically easier to provide explicit formulas for $P$ and $R$ rather 
than for $P$ and $Q$.

We start with the series $G$.


\subsection{Series $G$}\label{sec:G}

The polynomial $P$ for the series $G$ takes the form
\begin{eqnarray}\label{P-G}
P \eq A(x)^m,
\end{eqnarray}
where $A$ is a polynomial of degree $k-1$ whose roots are the black vertices (all of them 
are of degree~$m$). Notice that the number of these vertices does not coincide with the 
degree of the central vertex since we have one ``double'' edge.

We choose the normalization of $P$, $Q$ and $R=P-Q$ in the following way:
\begin{itemize}
\item	$P=A^m$ where $A$ is monic, $\deg A = k-1$;
\item	the central vertex is placed at $x=0$, so that $Q=x^k\cdot B$ where $B$ is  
		monic, $\deg B = n-k$; the roots of $B$ are the white vertices distinct 
		from zero;
\item	$R=c\,(x-1)$; this means that the pole inside the only face of degree 1 is 
		placed at $x=1$.
\end{itemize}

Thus, we get
\begin{equation} 
\label{xxyy}
A^m - c\,(x-1) \eq x^k\cdot B.
\end{equation} 

\begin{proposition}
The polynomial $A$ satisfies the differential equation 
\begin{eqnarray}\label{eq-for-AG} 
mA^{\prime}\cdot(x-1)-A=(m(k-1)-1)x^{k-1}.
\end{eqnarray} 
Consequently, coefficients $a_0,\ldots,a_{k-1}$ of $A(x) = \sum_{i=0}^{k-1}a_ix^i$
may be found by the following backward recurrence:
\begin{eqnarray}\label{AG}
a_{k-1} \eq 1, \qquad a_i \eq \frac{m(i+1)}{mi-1}\cdot a_{i+1} 
\qquad \mbox{for} \qquad 0\leq i\leq k-2. 
\end{eqnarray}
Finally, $c=-a_0$.
\end{proposition}

\paragraph{Proof.} Taking the derivative of the both sides of equality \eqref{xxyy} 
we obtain the equality 
$$
mA^{m-1}A^{\prime}-c \eq x^{k-1}\left(kB+xB^{\prime}\right),
$$ 
implying the equality
$$
mA^{m}A^{\prime}-cA \eq x^{k-1}A\left(kB+xB^{\prime}\right).
$$ 
Substituting in the last equality the value of $A^m$ from \eqref{xxyy}, we obtain 
$$
mA^{\prime}\left[c\,(x-1)+x^kB\right]-cA \eq x^{k-1}A\left(kB+xB^{\prime}\right)
$$
and 
$$
mA^{\prime}\cdot c\,(x-1)-cA \eq x^{k-1}\left[kAB+xAB^{\prime}-xmA^{\prime}B\right].
$$

We now observe that the degree of the left-hand side of the latter equality
is $k-1$, while its right-hand side is {\em proportional}\/ to $x^{k-1}$.
Therefore, the expression in the square brackets on the right is some constant $K$,
and both parts are equal to $K\cdot x^{k-1}$. The constant $K$ can be easily found 
as the leading coefficient of the left-hand side: it is equal to $mc(k-1)-c$.
Finally, we get the equality
$$
mcA^{\prime}-cA \eq (mc(k-1)-c)x^{k-1},
$$
which implies \eqref{eq-for-AG}. 

Substituting $A(x) = \sum_{i=0}^{k-1}a_ix^i$ in \eqref{eq-for-AG} we obtain 
\eqref{AG}. Finally, substituting $x=0$ in \eqref{xxyy} we obtain $c=-a_0^m$.
\hfill$\Box$

\begin{example}
Let us take $k=6$, so that $\deg A = k-1 = 5$. Then the corresponding polynomial 
looks as follows:
\begin{eqnarray}\label{A-F}
A = a_5x^5 + a_4x^4 + a_3x^3 + a_2x^2 + a_1x + a_0,
\end{eqnarray}
where
\begin{eqnarray*}\label{coeffs-G}
a_5 & = & 1, \\
a_4 & = & \frac{5m}{4m-1}\,, \\
a_3 & = & \frac{5m\cdot 4m}{(4m-1)(3m-1)}\,, \\
a_2 & = & \frac{5m\cdot 4m\cdot 3m}{(4m-1)(3m-1)(2m-1)}\,, \\
a_1 & = & \frac{5m\cdot 4m\cdot 3m\cdot 2m}{(4m-1)(3m-1)(2m-1)(m-1)}\,, \\
a_0 & = & \frac{5m\cdot 4m\cdot 3m\cdot 2m\cdot m}{(4m-1)(3m-1)(2m-1)(m-1)(-1)}\,
\end{eqnarray*}
\end{example}

\begin{remark}[Hypergeometric equation]\label{rem:hyp-geom}
Polynomial $A$ also satisfies the hypergeometric differential equation 
\begin{eqnarray}\label{eq2} 
x(1-x)\,\frac{d^2y}{dx^2} + \big[c-(a+b+1)x\big]\frac{dy}{dx} - ab\cdot y \eq 0.
\end{eqnarray}
Indeed, applying the differential operator 
$\displaystyle x\,\frac{d}{dx}+(1-k)$ to both parts of equality \eqref{eq-for-AG} 
we obtain
$$
x\left[mA^{\prime}\cdot(x-1)-A\right]^{\prime}+(1-k)\left[mA^{\prime}\cdot(x-1)-A\right]=0,
$$ 
implying 
$$
x(x-1)A^{\prime\prime} + \left[\left(1-\frac{1}{m}+(1-k)\right)x - 
(1-k)\right]A^{\prime} - \frac{(1-k)}{m}A \eq 0.
$$
Therefore, $A$ is a solution of the differential equation 
$$
x(1-x)\,\frac{d^2y}{dx^2} + \left[(1-k)-\left((1-k)-\frac{1}{m}+1\right)x\right]\frac{dy}{dx} 
+ \frac{(1-k)}{m}\,y \eq 0
$$
which is a particular case of \eqref{eq2} with 
$$
a=1-k, \quad b=-\frac{1}{m}, \quad c=1-k.
$$
\end{remark}

\subsection{Series $F$}\label{sec:F}

For this series we may assume that 
\begin{eqnarray}\label{P-F}
P \eq (x-1)^l A(x)^m, \qquad Q \eq x^kB(x).
\end{eqnarray}
Here $A$ is monic and $\deg A=k-1$; namely, $A$ is a polynomial whose roots are the 
black vertices of degree $m$. Now, $B$ is a polynomial whose roots are the white 
vertices distinct from zero, $\deg B=n-k$. The polynomials $P$ and $Q$ must 
satisfy the condition 
\begin{equation} \label{xxyyzz}
(x-1)^l A(x)^m-x^kB(x) \eq c\,, 
\end{equation}
where $c\in \C$ is a non-zero constant.

\begin{proposition}
The polynomial $A$ satisfies the differential equation 
\begin{eqnarray}\label{eq-for-AG2} 
mA^{\prime}\cdot(x-1)+lA \eq \left[m(k-1)+l\right]\,x^{k-1}.
\end{eqnarray} 
Consequently, coefficients $a_0,\ldots,a_{k-1}$ of $A(x) = \sum_{i=0}^{k-1}a_ix^i$
may be found by the following backward recurrence:
\begin{eqnarray}\label{AG0}
a_{k-1} \eq 1, \qquad a_i \eq \frac{m(i+1)}{mi+l}\cdot a_{i+1} 
\qquad \mbox{for} \qquad 0\leq i\leq k-2. 
\end{eqnarray}
Finally, the value of $c$ in \eqref{xxyyzz} is equal to $(-1)^la_0^m$.
\end{proposition}

\paragraph{Proof.} As above, let us take the derivative of both sides of 
equation \eqref{xxyyzz}. Then we get 
$$
(x-1)^{l-1}A^{m-1}\left[lA +m\,(x-1)A^{\prime}\right] \eq x^{k-1}\left(kB+xB^{\prime}\right).
$$ 
We observe that the polynomial $x^{k-1}$ is coprime with the factor $(x-1)^{l-1}A^{m-1}$
in the left-hand side, and therefore it must be proportional to the factor 
$lA +m\,(x-1)A^{\prime}$ which is itself a polynomial of degree $k-1$. Therefore,
both of them are equal to $K\cdot x^{k-1}$ where the constant $K$ can be found as
the leading coefficient of $lA +m\,(x-1)A^{\prime}$; namely, it is equal to
$m(k-1)+l$. Thus, \eqref{eq-for-AG2} holds.

Now, substituting $A(x) = \sum_{i=0}^{k-1}a_ix^i$ in \eqref{eq-for-AG2} we obtain
the recurrence \eqref{AG0}, and substituting $x=0$ in \eqref{xxyyzz} we obtain
the value of $c$.
\hfill$\Box$

\bsk

Here, like in the case of the series $G$, the polynomial $A$ also satisfies 
the hypergeometric differential equation, and therefore it may be represented 
through a hypergeometric function. 

\begin{example}
Let us take $k=6$, so that $\deg A = k-1 = 5$. Then the corresponding polynomial 
looks as follows:
\begin{eqnarray}
A = a_5x^5 + a_4x^4 + a_3x^3 + a_2x^2 + a_1x + a_0, \nonumber
\end{eqnarray}
where
\begin{eqnarray}\label{coeffs-F}
a_5 & = & 1\,, \nonumber \\
a_4 & = & \frac{5m}{l+4m}\,, \nonumber \\ 
a_3 & = & \frac{5m\cdot 4m}{(l+4m)(l+3m)}\,, \nonumber \\
a_2 & = & \frac{5m\cdot 4m\cdot 3m}{(l+4m)(l+3m)(l+2m)}\,, \nonumber \\
a_1 & = & \frac{5m\cdot 4m\cdot 3m\cdot 2m}{(l+4m)(l+3m)(l+2m)(l+m)}\,, \nonumber \\
a_0 & = & \frac{5m\cdot 4m\cdot 3m\cdot 2m\cdot m}{(l+4m)(l+3m)(l+2m)(l+m)l}\,. \nonumber
\end{eqnarray}
\end{example}


\subsection{Differential relations} 

The above method may be applied to DZ-pairs which do not necessary correspond to trees 
of diameter four or to unitrees. However, in general, it leads to {\it differential relations}\/
between $P$ and $Q$. Let us clarify what we mean by considering the problem of the 
difference between cubes and squares of polynomials, which was at the origin of the
whole activity concerning DZ-pairs, see \cite{BCHS-65}, \cite{Davenport-65}. 

Let $A,$  $B,$ and $R$  be polynomials such that
\begin{equation} \label{1} 
A^3-B^2 \eq R
\end{equation} 
and
$$
\deg A \eq 2k, \qquad \deg B \eq 3k, \qquad \deg R \eq k+1.
$$ 
Taking the derivative of both parts of \eqref{1} we obtain 
$$
3A^2A'-2BB' \eq R'.
$$ 
Multiplying now the last equality by $A$ and substituting $A^3$ from \eqref{1} we obtain 
the equality 
$$
3A'\left(B^2+R\right)-2BB'A \eq R'A,
$$ 
implying in its turn the equality 
$$
B\left(3A'B-2AB'\right) \eq R'A-3A'R.
$$
Since the degree of the right-hand side is
$$
\deg (R'A-3A'R)\leq 3k
$$ 
while $\deg B=3k$, the above equality implies that 
\begin{equation} \label{2} 
3A'B-2AB' \eq c
\end{equation}
for some non-zero constant $c\in\C$.

The last expression is a differential equation of the first order with respect to $A$ 
as well as with respect to $B$. Unfortunately, both $A$ and $B$ are unknown. Thus, 
it does not give us any immediate information about $A$ and $B$. Still, algebraic 
equations for coefficients of $A$ and $B$ obtained from \eqref{2} are (mostly) 
of degree~2 while the equations obtained from \eqref{1} are (mostly) of degree~3. 

Differentiating \eqref{2} and writing the expression thus obtained as a differential 
equation with respect to $A$ we get: 
\begin{equation} \label{3} 
A^{''} + \frac{B^{'}}{3B}\cdot A^{'} - \frac{\ \,2B^{''}}{3B}\cdot A \eq 0.
\end{equation} 
This differential equation is a particular case of the differential equation 
\begin{equation} \label{sti} 
\frac{d^2 S}{dz^2} + \left(\sum_{j=1}^m\frac{\gamma_j}{z-a_j}\right)\frac{dS}{dz}+\frac{V(z)}{\prod_{j=1}^m(z-a_j)}S \eq 0,
\end{equation}
where $V$ is a polynomial of degree at most $m-2.$ Polynomial solutions of the last equation 
are called Stieltjes polynomials. The polynomials $V$ for which \eqref{sti} has a polynomial 
solution are called Van Vleck polynomials. Thus, $B$ is a Van Vleck polynomial, and $A$ is the 
corresponding Stieltjes polynomial. 

Writing now \eqref{3} in the form 
$$
B^{''} - \frac{A^{'}}{2A}\cdot B^{'} - \frac{\ \, 3A^{''}}{2A}\cdot B \eq 0
$$ 
we obtain that $A$ is a Van Vleck polynomial and $B$ is the corresponding 
Stieltjes polynomial.

The above observations show that the relations between DZ-pairs and differential equations 
may be deeper than it seems at first glance and deserve further investigation.


\section{Series $H$ and $I$: decomposable ordinary trees}\label{sec:H}

\begin{figure}[htbp]
\begin{center}
\epsfig{file=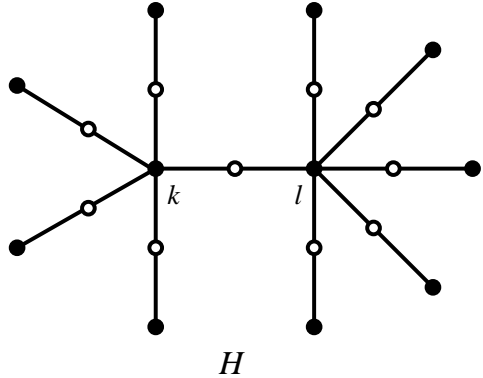,width=6.3cm}
\end{center}
\caption{Series $H$: ordinary trees of diameter 6 which are decomposable.}
\label{fig:H}
\end{figure}

In this section we consider series $H$ (Fig.~\ref{fig:H}) and $I$ (Fig.~\ref{fig:I}). 
In both cases the corresponding DZ-pairs are obtained with the help of the operation 
of composition. Notice that the trees in question are ordinary (the weights of all
edges are equal to 1). As it was mentioned in Definition~\ref{def:weighted}, Belyi
functions for ordinary trees are called Shabat polynomials.


\subsection{Series $H$}

The trees of the series $H$ are compositions of trees  from the series~$C$ with the 
parameters $s=t=1$ and chains of length~2.

\begin{figure}[t]
\begin{center}
\epsfig{file=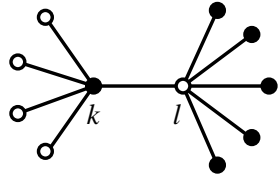,width=3.6cm}
\end{center}
\caption{Replace every edge of this tree with a two-edge chain, and you get the tree~$H$}
\label{fig:C-H}
\end{figure}

The expressions of the Shabat polynomials for the trees  from the series~$C$ in terms of
Jacobi polynomials are given in Sect.~\ref{sec:C}. Using the fact that $s=t=1$
we can also compute them directly. Indeed, the trees in question have exactly two 
vertices of degree greater than~1. Putting them into the points $x=0$ and $x=1$ and 
taking into account that the degree of the corresponding Shabat polynomial 
$S(x)$ is $k+l-1$, we conclude that the derivative of $S$ is proportional to 
$x^{k-1}(1-x)^{l-1}$. Therefore, the polynomial $S(x)$ itself can be written as
\begin{eqnarray}\label{S-H}
S(x) \eq K \cdot \int_0^x t^{k-1}(1-t)^{l-1}dt.
\end{eqnarray}
Then we automatically have $S(0)=0$, while in order to get $S(1)=1$ we must take
\begin{eqnarray}\label{K-H}
K \eq \frac{1}{B(k,l)}=\frac{(k+l-1)!}{(k-1)!(l-1)!}\,,
\end{eqnarray}
where
\begin{eqnarray}\label{euler-beta}
B(k,l) \eq \int_0^1 t^{k-1}(1-t)^{l-1}dt
\end{eqnarray}
is the Euler beta function.

Then, taking the Shabat polynomial for the chain with two edges and with two black
vertices put to 0 and 1, which is equal to
\begin{eqnarray}\label{U-H}
U(y) \eq 4y(1-y),
\end{eqnarray}
we obtain the following

\begin{proposition}
The polynomial $P$ for the tree $H$ is equal to
\begin{eqnarray}\label{P-H}
P(x) \eq U(S(x))
\end{eqnarray}
where $U$ is as in \eqref{U-H} and $S$ is as in \eqref{S-H} and \eqref{K-H}.
\end{proposition}

The proof is obvious.


\subsection{Series $I$}

\begin{figure}[htbp]
\begin{center}
\epsfig{file=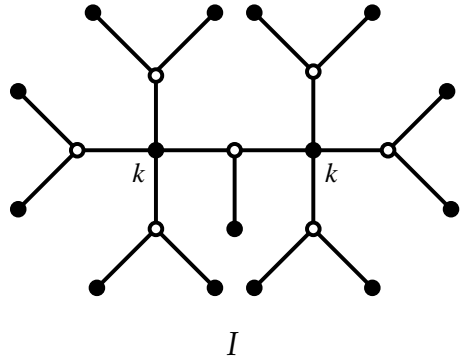,width=6cm}
\end{center}
\caption{Series $I$}\label{fig:I}
\end{figure}

Below are given Shabat polynomials $P(z)$ for the trees of the series $I$. These trees are
compositions of trees from the series~$C$  with $s=t=1$ and $k=l$, and the stars with 
three edges. Thus, $P(x) = U(S(x))$, where $S$ is a Shabat polynomial corresponding 
to a tree from the series $C$, and $U$ is a Shabat polynomial corresponding to the star 
with three edges. However, in order to achieve the rationality of the coefficients of 
$P$ we still must find an appropriate normalization of $S$.

\begin{figure}[htbp]
\begin{center}
\epsfig{file=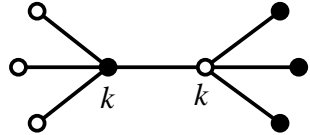,width=4cm}
\end{center}
\caption{Replace every edge with a three-edge star, and you get the tree~$I$}\label{fig:3kk-I}
\end{figure}

For this purpose, contrary to all traditions, let us put the vertices of degree~$k$ of 
the tree from the series~$C$ into the points $x=\pm\sqrt{-3}$. Then the derivative 
of the corresponding Shabat polynomials $S(x)$ must be equal to
\begin{eqnarray}\label{dif-S-I}
S'(x) \eq a\,(x+\sqrt{-3})^{k-1}(x-\sqrt{-3})^{k-1} \eq a\,(x^2+3)^{k-1}, \qquad a\in \C.
\end{eqnarray}
Therefore, 
\begin{eqnarray}\label{T-I}
S(x) \eq a\int (x^2+3)^{k-1}dx+b = a\left[\sum_{i=0}^{k-1}
\binom{k-1}{i}\frac{x^{2i+1}}{2i+1}3^{k-1-i}\right]+b
\end{eqnarray}
for some $b\in \C$. Substituting into $S(x)$ the critical
points $x=\pm\sqrt{-3}$, we obtain the critical values $b\pm c\sqrt{-3}$,
where 
\begin{equation} \label{x} 
c \eq a\cdot 3^{k-1}\sum_{i=0}^{k-1}\binom{k-1}{i}\frac{(-1)^{i}}{2i+1}.
\end{equation}

Setting  
\begin{equation} \label{y} 
b \eq -\frac{1}{2}
\end{equation} 
and choosing $a$ in such a way that 
\begin{equation} \label{z} 
c \eq \frac{1}{2}\,,
\end{equation} 
we obtain a polynomial $S\in\Q[x]$ with two critical values
\begin{eqnarray}\label{val-I}
y_{1,2} \eq \frac{-1\pm\sqrt{-3}}{2}.
\end{eqnarray}
Taking now
\begin{eqnarray}\label{U-I}
U(y) \eq 1-y^3
\end{eqnarray}
(we must take $1-y^3$ instead of $y^3$ in order to get the colors of the vertices
which would correspond to Fig.~\ref{fig:I}), we obtain the following

\begin{proposition}
The polynomial $P(x)$ for the tree $I$ is equal to
\begin{eqnarray}\label{P-I}
P(x) \eq U(S(x)),
\end{eqnarray}
where $U$ is as in \eqref{U-I} and $S$ is as in \eqref{T-I} with $a$ and $b$ defined
by conditions \eqref{x}, \eqref{y}, \eqref{z}.
\end{proposition}

Once again, the proof is obvious.


\section{Series $J$}\label{sec:J}

\begin{figure}[htbp]
\begin{center}
\epsfig{file=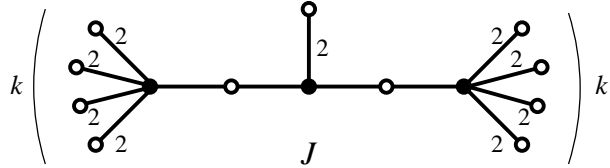,width=8cm}
\end{center}
\caption{Series $J$}\label{fig:J}
\end{figure}

This is the last infinite series of unitrees. The degree of this tree, or its
total weight, is $2k+6$.

Let us normalize the polynomial $P$ so that 
\begin{eqnarray}
P \eq (x+1)^4\cdot (x^2+a)^{2k+1}.
\end{eqnarray}
This means that the black vertex of degree 4 is put at $x=-1$, while two black
vertices of degree $2k+1$ are put at the points $\pm\sqrt{-a}$ for certain $a\in\Q$,
$a>0$.

All the white vertices are of degree 2; therefore, the polynomial $Q$ has the form 
$$
Q(x) \eq A(x)^2
$$ 
for some polynomial $A$, $\deg A=2k+3$. Further, condition \eqref{dere} gives us 
$$
(x+1)^4\cdot (x^2+a)^{2k+1}-A(x)^2\, \rav_{x\to \infty}\,O(x^{2k+1});
$$
here $2k+1$ is the ``overweight'' of the tree (that is, its total weight minus the
number of edges of the topological tree). For the reciprocal polynomials
this gives (see \eqref{dere0})
\begin{equation} \label{kon} 
P^*-Q^* \eq (1+x)^4\cdot (1+ax^2)^{2k+1}-A^{*}(x)^2\,\rav_{x\to 0}\,O(x^{2k+5});
\end{equation} 
here $2k+5$ is the number of edges of the topological tree.

\begin{proposition}
The reciprocal polynomials $P^*$ and $Q^*$ may be represented as follows:
\begin{eqnarray}\label{P*}
P^* \eq (1+x)^4\cdot (1+(2k+4)x^2)^{2k+1}, \qquad Q^*(x) \eq A^{*}(x)^2,
\end{eqnarray}
where $A^*$  is the initial segment of the series $(P^*)^{1/2}$ up to the degree $2k+3$:
\begin{eqnarray}
(1+x)^2(1+(2k+4)x^2)^{(2k+1)/2}\,\rav_{x\to 0}\, A^* + O(x^{2k+4}).
\end{eqnarray}
\end{proposition}

\paragraph{Proof.} 
Let \begin{eqnarray}\label{A*}
A^* \eq (1+x)^2\cdot (1+ax^2)^{(2k+1)/2} + x^{2k+4}\cdot h
\end{eqnarray}
where 
$$
h\,\rav_{x\to 0}\,O(1).
$$ 
Computing $(A^*)^2$ we get
{\small
\begin{eqnarray}\label{raznost}
(A^*)^2 & = & P^* + 2x^{2k+4}\cdot h\cdot (1+x)^2\cdot (1+ax^2)^{(2k+1)/2} + 
              x^{4k+8}\cdot h^2 \nonumber \\
        & = & P^* + x^{2k+4}\left[2h\cdot (1+x)^2\cdot (1+ax^2)^{(2k+1)/2} + 
		      x^{2k+4}\cdot h^2\right].
\end{eqnarray}
}

\ni
Thus, for any value of the parameter $a$ we have
$$
P^*-A^{*}(x)^2\, \rav_{x\to 0}\,O(x^{2x+4}),
$$ 
and therefore, in order to obtain \eqref{kon}, we only have to show that for $a=2k+4$ 
the constant term of $h$ is equal to zero, or, equivalently, the coefficient 
in front of $x^{2k+4}$ in the series 
$$
(P^*)^{1/2} \eq (1+x)^2\cdot (1+ax^2)^{(2k+1)/2}
$$ 
vanishes.

Let us write the second factor of the latter expression explicitly:
{\small
\begin{eqnarray}\label{binom}
(1+ax^2)^{(2k+1)/2} & = & 1 + \frac{2k+1}{2}ax^2 + 
\frac{1}{2!}\cdot\frac{(2k+1)(2k-1)}{4}a^2x^4 + \nonumber \\
& & \frac{1}{3!}\cdot\frac{(2k+1)(2k-1)(2k-3)}{8}a^3x^6 + \ldots + \nonumber \\
& & \frac{1}{(k+2)!}\frac{(2k+1)(2k-1)\ldots(-1)}{2^{k+2}}a^{k+2}x^{2k+4} + \ldots 
\end{eqnarray}}
Notice that this series involves only even powers. Multiplying it by 
$$
(1+x)^2 \eq 1+2x+x^2
$$ 
we see that the coefficient in front of $x^{2k+4}$ in $(P^*)^{1/2}$
is the sum of the coefficients in front of $x^{2k+4}$ and $x^{2k+2}$ in \eqref{binom}. 
Therefore, we must ensure that
\begin{eqnarray}\label{coeff:2k+4}
\frac{1}{(k+1)!}\cdot\frac{(2k+1)(2k-1)\ldots\cdot 1}{2^{k+1}}\cdot a^{k+1}\,\, +
& & \nonumber \\
\frac{1}{(k+2)!}\cdot\frac{(2k+1)(2k-1)\ldots\cdot(-1)}{2^{k+2}}\cdot a^{k+2} 
& = &  0.
\end{eqnarray}
Collecting similar terms we get
\begin{eqnarray}
\frac{1}{(k+1)!}\cdot\frac{(2k+1)(2k-1)\ldots\cdot 1}{2^{k+1}}\cdot a^{k+1}
\left(1+\frac{1}{k+2}\cdot\frac{(-1)}{2}\cdot a \right) \eq 0,
\end{eqnarray}
which gives $a=2k+4$.
\hfill$\Box$

\begin{example}
Let us take $k=3$. Then we have:
\begin{eqnarray}\label{P3-J}
P^* \eq (1+x)^4 (1+10x^2)^7.\nonumber
\end{eqnarray}
Further, 
\begin{eqnarray}\label{taylor-J}
(P^*)^{1/2} & = & (1+x)^2 (1+10\,x^2)^{7/2} \nonumber \\ 
            & = & 1 + 2\,x + 36\,x^2 + 70\,x^3 + \frac{945}{2}\,x^4 + 875\,x^5 + 2625\,x^6 +
                  4375\,x^7 + \nonumber \\
            &   & \frac{39\,375}{8}\,x^8 + \frac{21\,875}{4}\,x^9 - 
                  \frac{21\,875}{4}\,x^{11} + \frac{65\,625}{16}\,x^{12} + \ldots \nonumber
\end{eqnarray}
Notice that the term with $x^{10}$ is missing. Finally,
\begin{eqnarray}
A^* & = & 1 + 2\,x + 36\,x^2 + 70\,x^3 + \frac{945}{2} \,x^4 + 875\,x^5 + 2625\,x^6 +
4375\,x^7 + \nonumber \\
& & \frac{39\,375}{8}\,x^8 + \frac{21\,875}{4}\,x^9. \nonumber
\end{eqnarray}

\end{example}


\section{Sporadic trees}\label{sec:sporadic}

As it was explained previously, in Sect.~\ref{sec:computation}, the verification of
the results given below is trivial. Therefore, we present nothing else but
the polynomials themselves.

\subsection{Tree $K$}

\begin{figure}[htbp]
\begin{center}
\epsfig{file=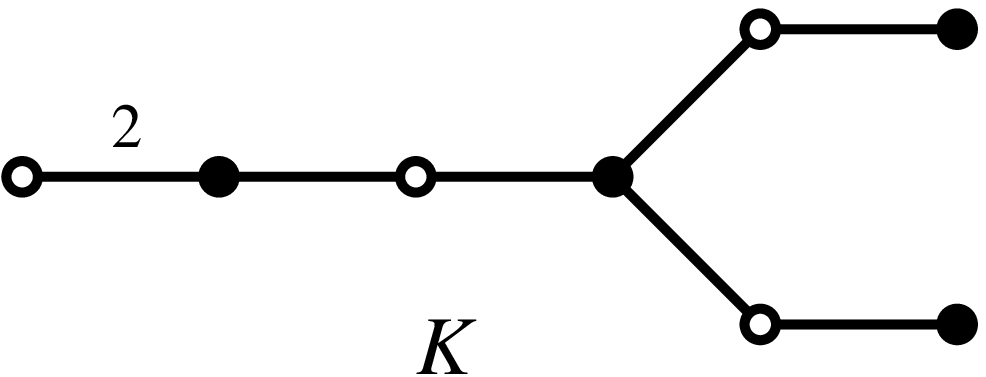,width=5cm}
\end{center}
\caption{Tree $K$}\label{fig:K}
\end{figure}

\begin{eqnarray}
P & = & (x^2-5x+1)^3(x^2-13x+49), \nonumber \\
Q & = & (x^4-14x^3+63x^2-70x-7)^2, \nonumber \\
R & = & -1728x. \nonumber
\end{eqnarray}

\subsection{Tree $L$}

\begin{figure}[htbp]
\begin{center}
\epsfig{file=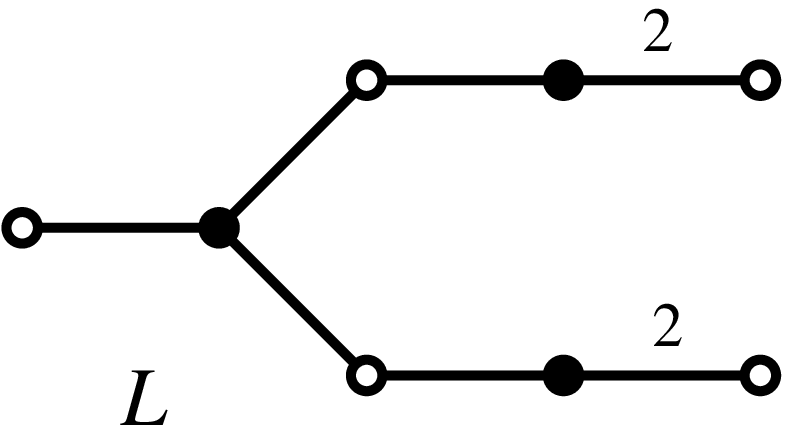,width=5cm}
\end{center}
\caption{Tree $L$}\label{fig:L}
\end{figure}

\begin{eqnarray}
P & = & (x^3-16x^2+160x-384)^3, \nonumber \\
Q & = & x\,(x^4-24x^3+336x^2-2240x+8064)^2, \nonumber \\
R & = & -2^{14}\cdot 3^3\,(x^2-13x+128). \nonumber
\end{eqnarray}

\subsection{Tree $M$}

\begin{figure}[htbp]
\begin{center}
\epsfig{file=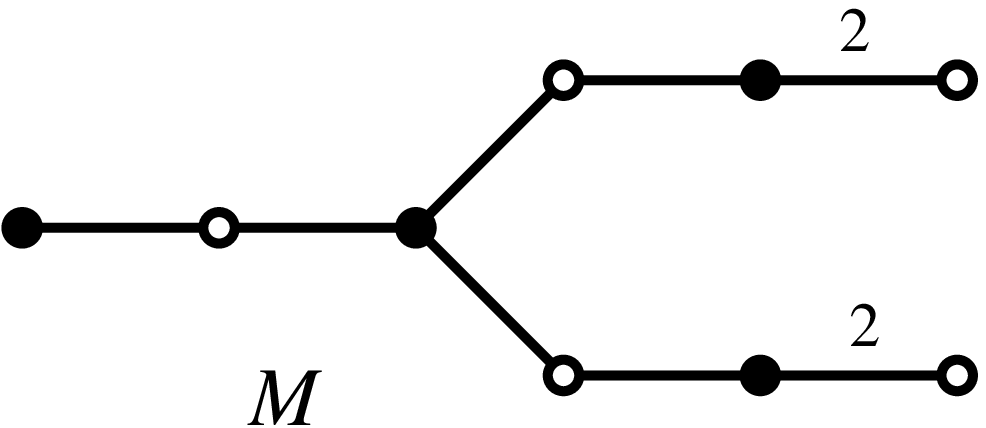,width=6cm}
\end{center}
\caption{Tree $M$}\label{fig:M}
\end{figure}

\begin{eqnarray}
P & = & x\,(x^3-36x^2+540x-2592)^3, \nonumber \\
Q & = & (x^5-54x^4+1296x^3-15\,552x^2+87\,480x+104\,976)^2, \nonumber \\
R & = & -2^6\cdot 3^{12}\,(x^2-28x+324). \nonumber
\end{eqnarray}

\subsection{Tree $N$}

\begin{figure}[htbp]
\begin{center}
\epsfig{file=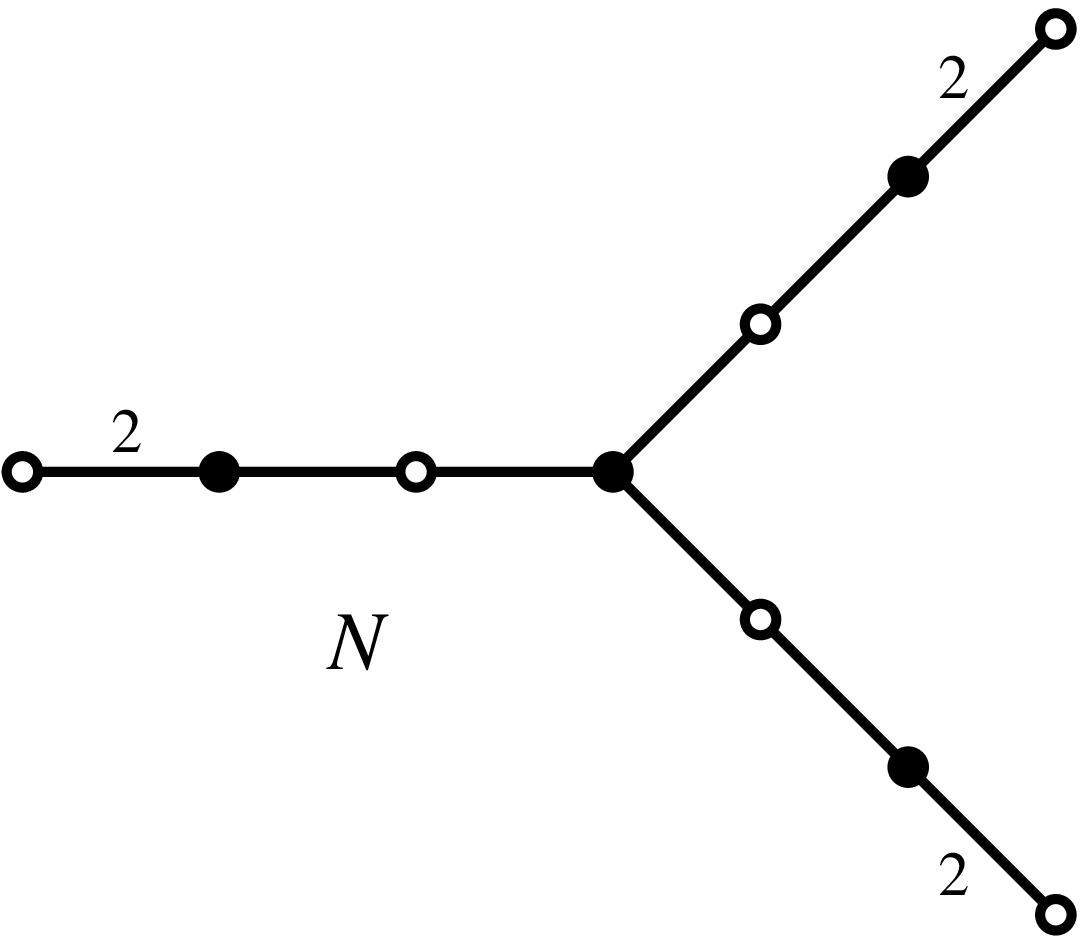,width=6cm}
\end{center}
\caption{Tree $N$}\label{fig:N}
\end{figure}

\begin{eqnarray}
P & = & x^3\,(x^3-8)^3, \nonumber \\
Q & = & (x^6-12x^3+24)^2, \nonumber \\
R & = & 64\,(x^3-9). \nonumber
\end{eqnarray}

\ni
This tree is symmetric, with the symmetry of order~3. Therefore, $P$, $Q$, $R$ are 
polynomials in~$x^3$.

\subsection{Tree $O$}\label{sec:O}

\begin{figure}[htbp]
\begin{center}
\epsfig{file=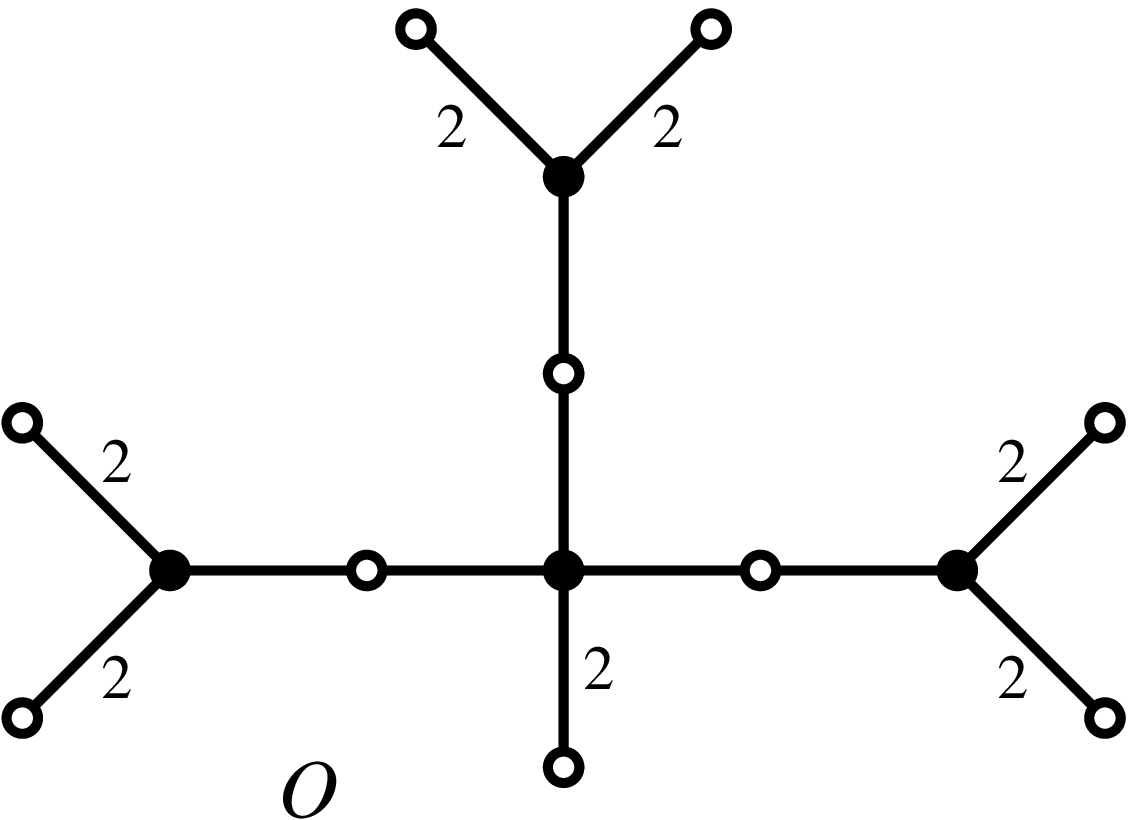,width=6.5cm}
\end{center}
\caption{Tree $O$}\label{fig:O}
\end{figure}

\begin{eqnarray}
P & = & (x^4+6x^2+64x-55)^5, \nonumber \\
Q & = & (x^{10}+15x^8+160x^7-70x^6+1440x^5+6510x^4 \nonumber \\
  &   & -\, 11\,040x^3+26\,805x^2+40\,160x-226\,797)^2, \nonumber \\
R & = & 2^{20}\,(5x^7+59x^5+690x^4-485x^3+3820x^2 \nonumber \\
  &   & +\,20\,165x-49\,534). \nonumber
\end{eqnarray}

This triple was found in Beukers and Stewart \cite{BeuSte-10} (only the polynomial~$P$
is given in their paper, but it uniquely determines two other polynomials).

\subsection{Tree $P$}\label{sec:sporadic-P}

\begin{figure}[htbp]
\begin{center}
\epsfig{file=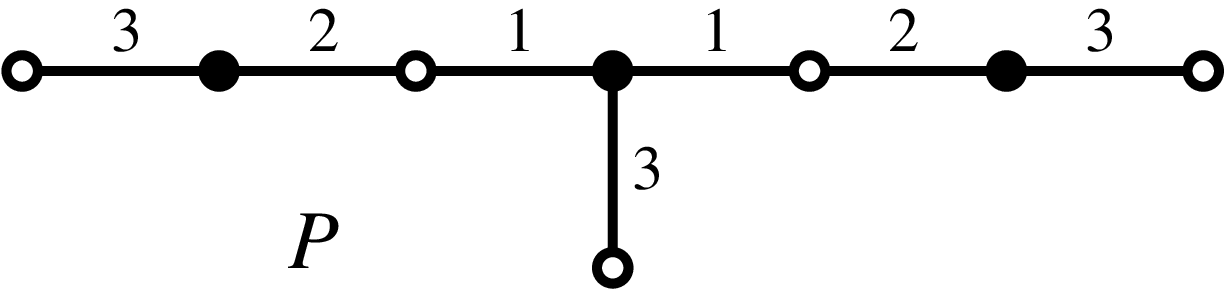,width=7cm}
\end{center}
\caption{Tree $P$}\label{fig:P}
\end{figure}

\begin{eqnarray}
P & = & (x^3+9x+9)^5, \nonumber \\
Q & = & (x^5+15x^3+15x^2+45x+90)^3, \nonumber \\
R & = & -27\,(15x^8+395x^6+423x^5+3330x^4+7290x^3 \nonumber \\
  &   & +\,11\,880x^2+29\,565x+24\,813). \nonumber
\end{eqnarray}

Once again, the answer is taken from \cite{BeuSte-10}, with a slight renormalization.

\subsection{Tree $Q$}

\begin{figure}[htbp]
\begin{center}
\epsfig{file=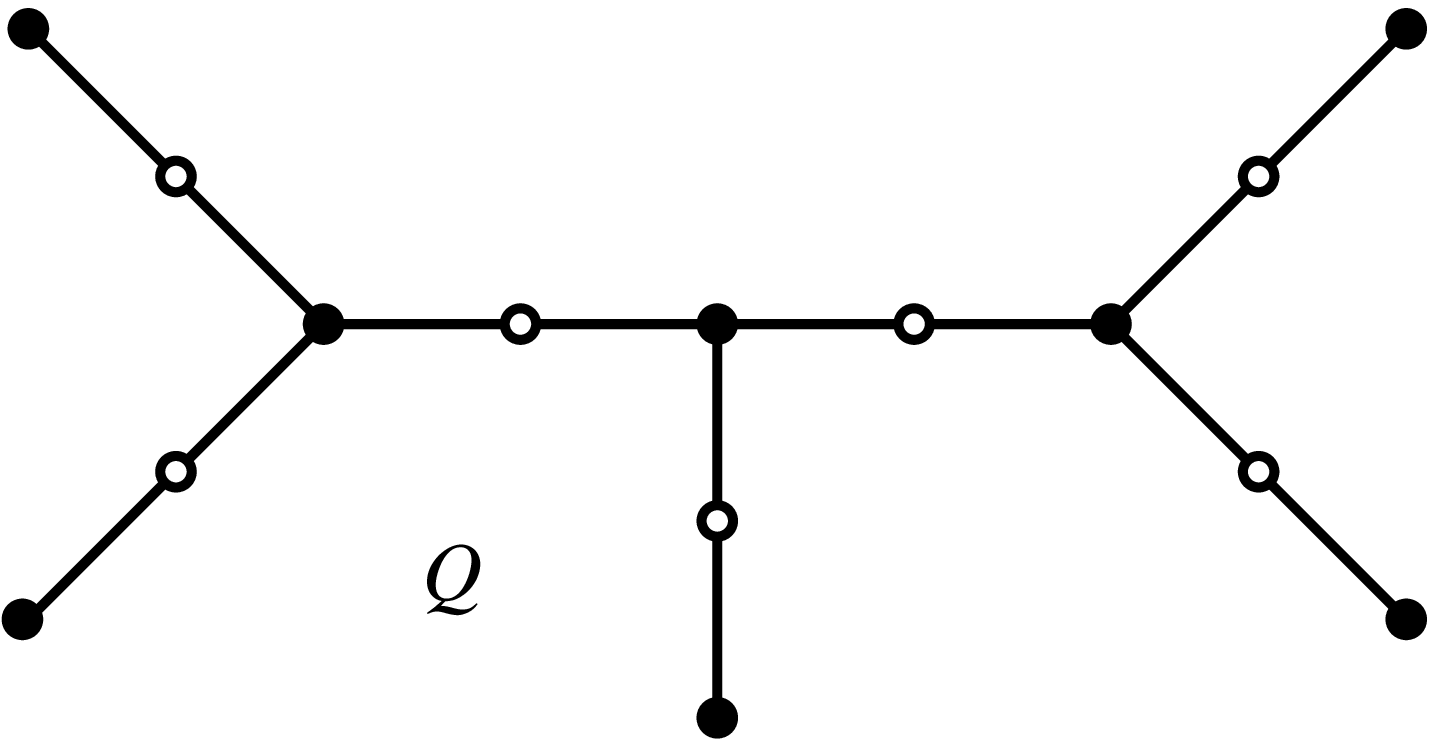,width=8cm}
\end{center}
\caption{Tree $Q$}\label{fig:Q}
\end{figure}

\begin{eqnarray}
P & = & (x^3+15x+16)^3 (x^5+39x^3+64x^2+384x+1872), \nonumber \\
Q & = & (x^7+42x^5+56x^4+525x^3+1680x^2+1792x+6456)^2, \nonumber \\
R & = & -2^6\cdot 3^{12}. \nonumber
\end{eqnarray}

This tree is the only sporadic tree from the Adrianov's list of ordinary unitrees. 
Correspondingly, $P$ is a Shabat polynomial: the polynomial $R$ is a constant.

\bigskip

Note that the positions of certain black vertices are rational:
$$
x^3+15x+16 = (x+1)(x^2-x+16),
$$
$$
x^5+39x^3+64x^2+384x+1872 = (x+3)(x^4-3x^3+48x^2-80x+624).
$$

\subsection{Tree $R$}

\begin{figure}[htbp]
\begin{center}
\epsfig{file=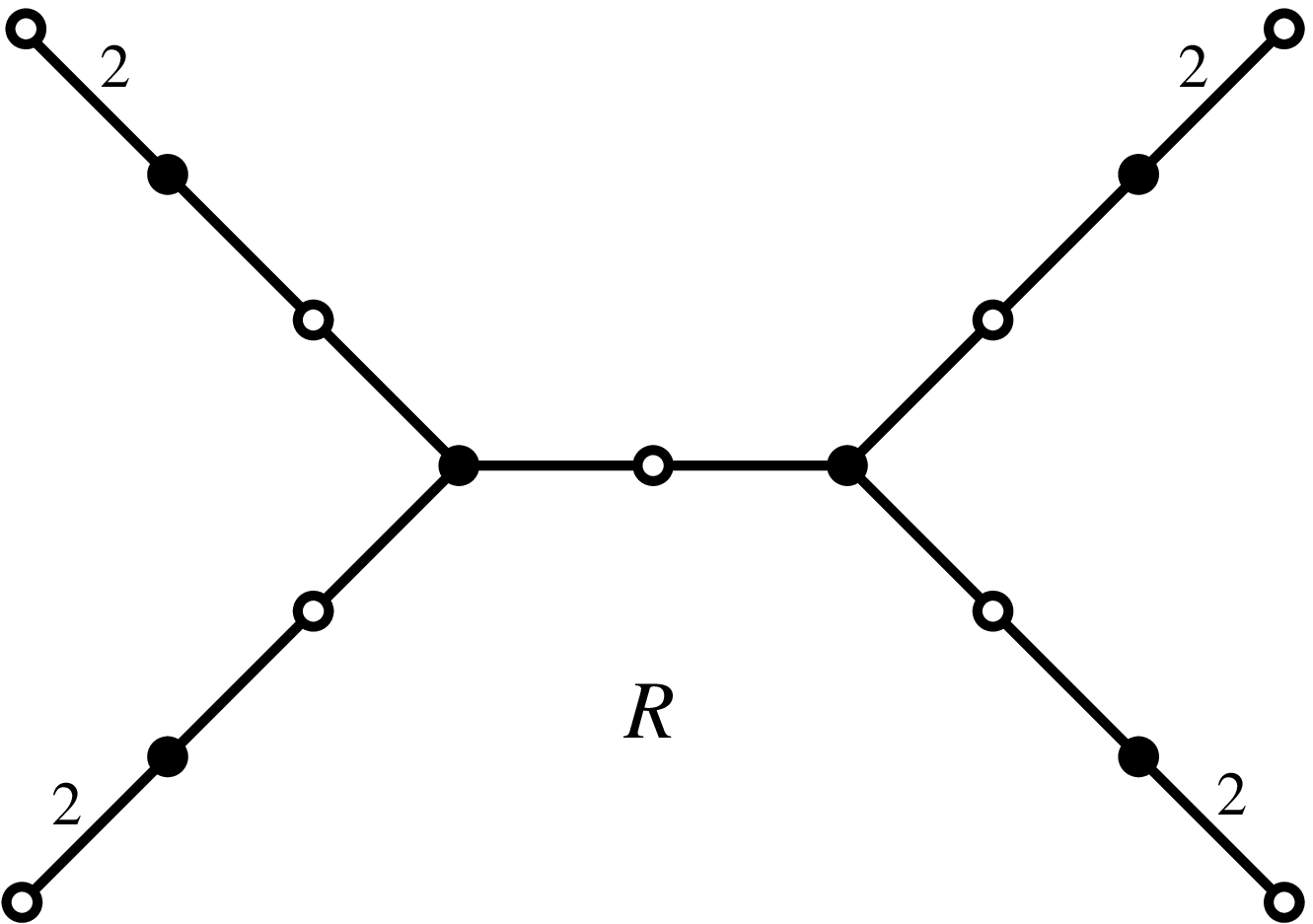,width=8cm}
\end{center}
\caption{Tree $R$}\label{fig:R}
\end{figure}

The tree $R$ is the ``square'' of the tree $L$: it is symmetric, with the symmetry
of order~2, and one of its ``halves'' is equal to $L$. Therefore, we may take the 
polynomials for the tree $L$ and insert $x^2$ instead of $x$.

\begin{eqnarray}
P & = & (x^6-16x^4+160x^2-384)^3, \nonumber \\
Q & = & x^2\,(x^8-24x^6+336x^4-2240x^2+8064)^2, \nonumber \\
R & = & -2^{14}\cdot 3^3\,(x^4-13x^2+128). \nonumber
\end{eqnarray}

\subsection{Tree $S$}\label{sec:tree-S}

\begin{figure}[htbp]
\begin{center}
\epsfig{file=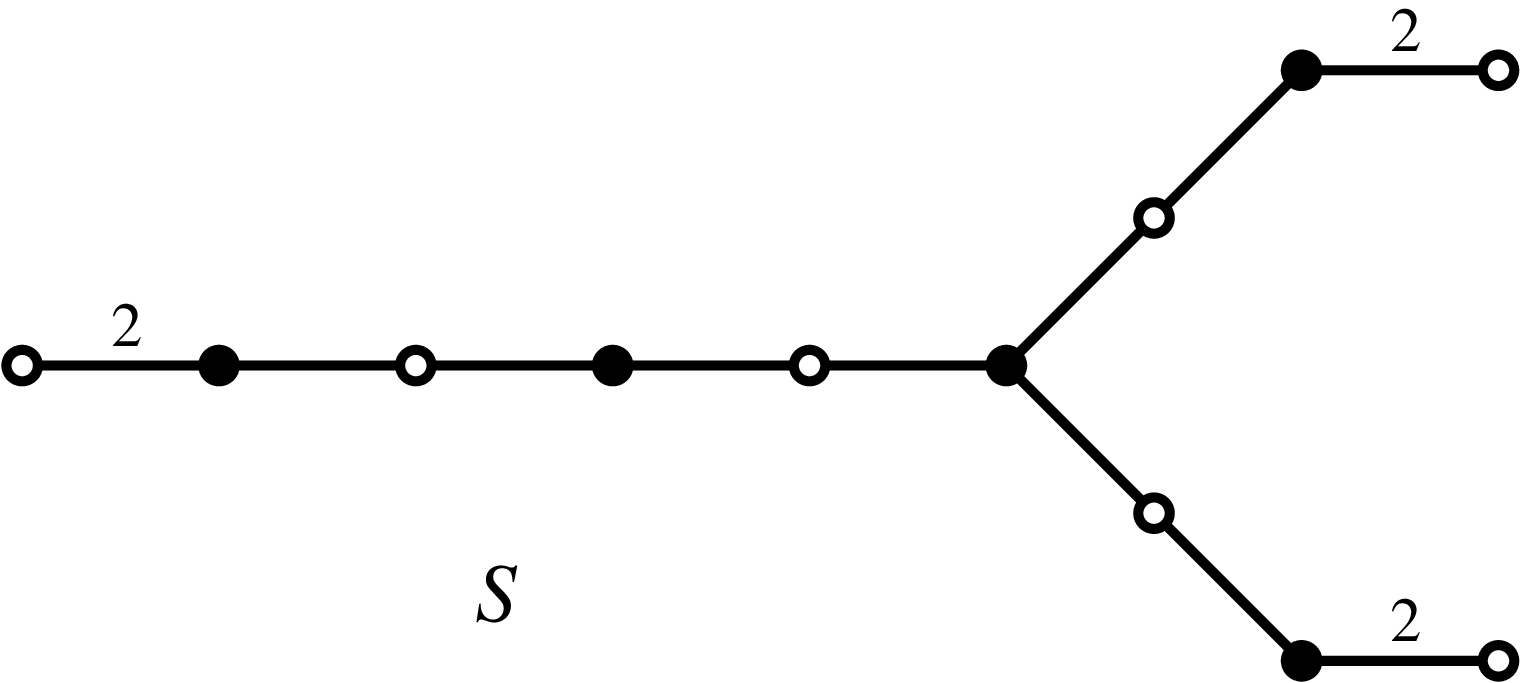,width=8cm}
\end{center}
\caption{Tree $S$}\label{fig:S}
\end{figure}

\begin{eqnarray}
P & = & x^2\,(x^4+24x^3+176x^2-2816)^3, \nonumber \\
Q & = & (x^7+36x^6+480x^5+2304x^4-3840x^3, \nonumber \\
  &   &   -\,\,55\,296x^2-14\,336x+221\,184)^2 \nonumber \\
R & = & 2^{22}\cdot 3^3\,(x^3+17x^2+56x-432). \nonumber
\end{eqnarray}

Notice that the second factor in $P$, the one which is ``cubed'', does not contain the 
term with $x$: this is not a misprint.

\subsection{Tree $T$}

The picture of this tree is given in Example~\ref{ex:pair-T.bis}, and the corresponding polynomials 
are given in Example~\ref{ex:pair-T}.


\section{Trees defined over $\Q$ by virtue of Galois \\ theory}\label{sec:galois}

Recall that the passport of a (bicolored weighted plane) tree is a pair of partitions 
$\al,\be\vdash n$, where $n$ is the degree (or the total weight) of the tree, $\al$
represents the set of degrees of its black vertices, and $\be$ represents the set of
degrees of its white vertices.

\begin{definition}[Combinatorial orbit]\label{def:comb-orbit}
The set of weighted trees with the same passport is called {\em combinatorial
orbit}.
\end{definition}

Unitrees represent, in fact, combinatorial orbits consisting of a unique tree.

Usually, a DZ-pair corresponding to a tree is defined over a number field
whose degree is equal to the size of the combinatorial orbit to which this tree
belongs. This is why unitrees are always defined over $\Q$. There exist, however,
other Galois invariants which may split a combinatorial orbit into several
distinct Galois orbits. In this way we may obtain certain trees which are not unitrees
but which are still defined over $\Q$. In \cite{PakZvo.I-14} we gave several such
examples. Here we present the corresponding DZ-pairs.

\subsection{A tree with the monodromy group ${\rm PGL}_2(7)$}

\begin{figure}[htbp]
\begin{center}
\epsfig{file=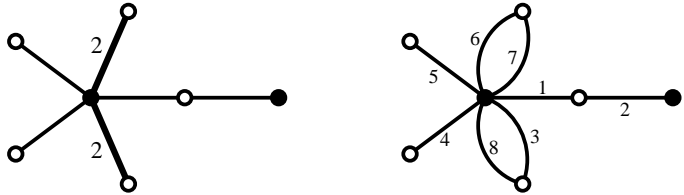,width=9cm}
\end{center}
\caption{The monodromy group of this tree is ${\rm PGL}_2(7)$. Numbers written on the
edges of the tree on the left are their weights; numbers written on the edges of the 
map on the right are {\em not}\/ weights: they are edge labels from 1 to 8.}
\label{fig:pgl27}
\end{figure}

A bicolored map may be characterized by a pair of permutations acting on the set of its
edges: one permutation represents the cyclic order (in the positive direction) of the edges 
around black vertices, the other one, the cyclic order around white vertices. For example,
the map shown in Fig.~\ref{fig:pgl27}, is represented by the pair of permutations
$$
a \eq (1,7,6,5,4,8,3), \qquad b \eq (1,2)(3,8)(6,7).
$$
It turns out that the permutation group $G=\langle a,b\rangle$ is equal to ${\rm PGL}_2(7)$.
This group, which is called {\em monodromy group}, is a Galois invariant. Since this tree
is the only one in its combinatorial orbit whose monodromy group is ${\rm PGL}_2(7)$, it is
defined over $\Q$.

\begin{eqnarray}
P & = & x^7(x-6), \nonumber \\
Q & = & (x^3-6x^2+12x-36)^2(x^2+6x+12), \nonumber \\
R & = & -2^4\cdot 3^3\,(7x^2-6x+36). \nonumber
\end{eqnarray}

The combinatorial orbit to which this tree belongs, that is, the set of trees
with the passport $(7^11^1,2^31^2)$, contains six trees. The five remaining 
trees constitute a single Galois orbit; the corresponding DZ-pairs (or, we
may say, the trees themselves) are defined over the splitting field of the 
polynomial
\begin{eqnarray}\label{def1}
a^5+22a^4+209a^3+1040a^2+2624a+2560. \nonumber
\end{eqnarray}

\subsection{Another tree with the monodromy group ${\rm PGL}_2(7)$}

\begin{figure}[htbp]
\begin{center}
\epsfig{file=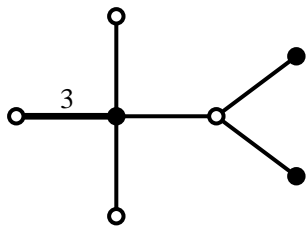,width=4cm}
\end{center}
\caption{This tree also has monodromy group ${\rm PGL}_2(7)$.}
\label{fig:pgl27-bis}
\end{figure}

The combinatorial orbit corresponding to the passport $(6^11^2,3^21^2)$,
consists of five trees. One of them, shown in Fig.~\ref{fig:pgl27-bis},
has the monodromy group ${\rm PGL}_2(7)$. Therefore, it is defined over $\Q$.
Its DZ-pair is given below.

\begin{eqnarray}
P & = & x^6(x^2-9x+21), \nonumber \\
Q & = & (x^2-3x-3)^3(x^2+3), \nonumber \\
R & = & 27\,(7x^2+9x+3). \nonumber
\end{eqnarray}

\begin{figure}[htbp]
\begin{center}
\epsfig{file=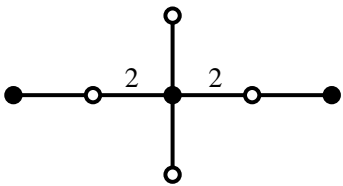,width=4.5cm}
\end{center}
\caption{The symmetric tree with the passport $(6^11^2,3^21^2,6^11^2)$.}
\label{fig:sym}
\end{figure}

One of the trees in this combinatorial orbit is symmetric (see Fig.~\ref{fig:sym}) 
and is therefore also defined over $\Q$.

The corresponding polynomials are
\begin{eqnarray}
P & = & x^6(x^2-2), \nonumber \\
Q & = & (x^2-1)^3(x^2+1), \nonumber \\
R & = & -2x^2+1. \nonumber
\end{eqnarray}

The three remaining trees constitute a single Galois orbit and are defined over the 
splitting field of the polynomial
\begin{eqnarray}\label{def2}
a^3-6a+16. \nonumber
\end{eqnarray}

\subsection{A series in which one of the trees is self-dual}\label{sec:self-dual}

The {\em duality}\/ for the bicolored maps is defined as follows:
\begin{itemize}
\item   a map and its dual share their white vertices;
\item   black vertices of each map correspond to the faces of the dual map;
\item   edges of the dual map connect the centers of the faces of the initial map to
        the white vertices which lie on the border of these faces.
\end{itemize}
See details and examples in \cite{PakZvo.I-14}. A map is {\em self-dual}\/ if it is
isomorphic to its dual. Self-duality is a Galois invariant. The maps corresponding to
weighted trees may well be self-dual.

\msk

\begin{figure}[htbp]
\begin{center}
\epsfig{file=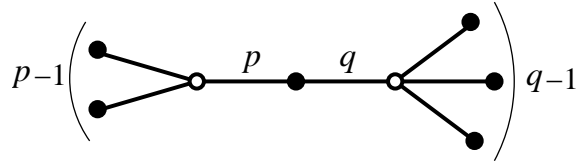,width=7.5cm}
\end{center}
\caption{Self-dual tree.}\label{fig:self-dual}
\end{figure}

Let $p,q$ be two positive integers, and $p<q$. We consider the trees with the black
partition $\al=(p+q,1^{p+q-2})$ and the white partition $\be=(2p-1,2q-1)$. The partition
representing the face degrees is $\ga=(p+q,1^{p+q-2})$. We notice that $\ga=\al$;
therefore, the corresponding combinatorial orbit {\em may}\/ contain self-dual
trees. It is easy to verify that this combinatorial orbit consists of $2p-1$ trees, 
and that only one of them is self-dual, namely, the tree shown in Fig.~\ref{fig:self-dual}. 
Therefore, this tree is defined over $\Q$.

Put the white vertices at the points $x=-1$ and $x=1$ so that 
\begin{eqnarray}\label{Q4}
Q(x) \eq (x+1)^{2p-1}\,(x-1)^{2q-1} \nonumber
\end{eqnarray}
(notice that both powers are odd). Observe now that this polynomial is ``antipalindromic'': 
if we write it as
$$
a_nx^n+a_{n-1}x^{n-1}+\ldots+a_1x+a_0, 
$$
then $a_n = -a_0$, $a_{n-1} = -a_1$, \ldots\ 
This fact trivially follows from the equality $x^n\cdot Q(1/x)=-Q(x)$. Because of this,
the coefficient in front of the ``middle'' degree $n/2=p+q-1$ is zero. Therefore,
if we take the higher degrees from $2p+2q-2$ to $p+q$, what will remain is a polynomial
of degree $p+q-2$. In other words,
\begin{eqnarray}\label{poly}
Q(x) \eq x^{p+q}\cdot A(x) - R(x), \nonumber
\end{eqnarray}
where $\deg A=\deg R = p+q-2$.
Setting now 
\begin{eqnarray}\label{P4}
P(x) \eq x^{p+q}\cdot A(x), \nonumber
\end{eqnarray} 
we see that $P$, $Q$ is a DZ-pair with required properties.
Notice that the polynomial $R(x)$ is reciprocal to $A(x)$. Geometrically, this means that
if $x_1,x_2,\ldots,x_m$ are the positions of the black vertices of degree~1 (here
$m=p+q-2$), then the centers of the faces of degree~1 are $1/x_1,1/x_2,\ldots,1/x_m$.
Together with the fact that the position of the black vertex of degree $p+q$ is $x=0$
while the center of the face of degree $p+q$ is $\infty$, this shows that the
map in question is indeed self-dual.

\begin{example}
Let us take, for example, $p=2$, $q=5$. Then
\begin{eqnarray}
Q(x)=(x+1)^3\,(x-1)^9 & = & x^{12}-6x^{11}+12x^{10}-2x^9-27x^8+36x^7 \nonumber \\
                      &   & -\,\,(1-6x+12x^2-2x^3-27x^4+36x^5) \nonumber \\
                      & = & x^7\cdot A(x) - R(x) \eq P(x)-R(x), \nonumber
\end{eqnarray}
where $\deg A = \deg R = 5$ and $R=A^*$. 
\end{example}

\subsection{A ``historical'' sporadic example}

\begin{figure}[htbp]
\begin{center}
\epsfig{file=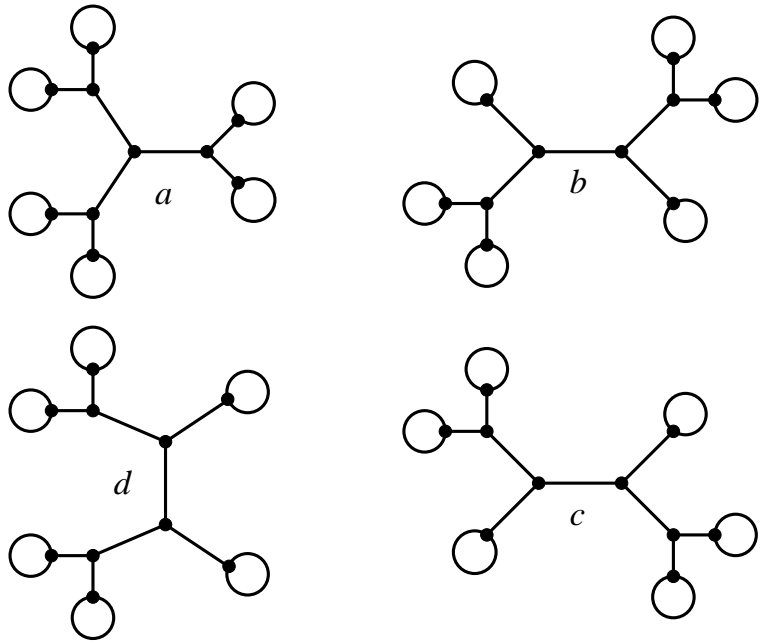,width=10cm}
\end{center}
\caption{Four trees with the passport $(3^{10},2^{15})$. The trees $a$ and $d$ are 
defined over $\Q$.}
\label{fig:zannier-5}
\end{figure}

The combinatorial orbit corresponding to the passport $\al=3^{10}$, $\be=2^{15}$
is shown in Fig.~\ref{fig:zannier-5}. It consists of four trees (recall that the
invisible white vertices are middle points of the edges), and is divided into three 
Galois orbits. 

The tree $a$ is the only one which is symmetric with the symmetry of order~3.
Therefore, it is defined over $\Q$. The corresponding polynomials were computed 
by B.\,Birch in 1965 \cite{BCHS-65}. They look as follows (notice that they are
polynomials in~$x^3$):
\begin{eqnarray}
P_a(x) & = & x^3 (x^9 + 12x^6 + 60x^3 + 96)^3\,, \nonumber \\
Q_a(x) & = & (x^{15} + 18x^{12} + 144x^9 + 576x^6 + 1080x^3 + 432)^2\,, \nonumber \\
R_a(x) & = & -1728\,(3x^6 + 28x^3 + 108)\,. \nonumber
\end{eqnarray}

The trees $b$ and $c$ are symmetric with the symmetry of order 2 with respect
to an (invisible) white vertex. They are also mirror
symmetric to each other; therefore, the complex conjugation sends one of the trees to
the other. Thus, we may conclude that this couple of trees constitutes a separate
Galois orbit, and this orbit is defined over an imaginary quadratic field. 
The corresponding polynomials were computed in 2005 by Shioda \cite{Shioda-05} and,
indeed, they are defined over the field $\Q(\sqrt{-3})$. We do not present these
polynomials here.

The tree $d$ does not have any particular combinatorial properties. (It is known
that the mirror symmetry of a dessin is not a Galois invariant.) But {\em it remains alone}, 
that is, it constitutes a Galois orbit containing a single element. Therefore,
it is defined over $\Q$. The corresponding polynomials were computed in 2000 by 
N.\,Elkies \cite{Elkies-00}. They look as follows:
\begin{eqnarray}
P_d(x) & = & (x^{10} - 2x^9 + 33x^8 - 12x^7 + 378x^6 + 336x^5 + 2862x^4 \nonumber \\ 
     &   &  + \,\, 2652x^3 + 14\,397x^2 + 9922x + 18\,553)^3, \nonumber \\
Q_d(x) & = & (x^{15} - 3x^{14} + 51x^{13} - 67x^{12} + 969x^{11} + 33x^{10} + 10\,963x^9 \nonumber \\
     &   &  +\,\, 9729x^8 + 96\,507x^7 + 108\,631x^6 + 580\,785x^5 + 700\,503x^4 \nonumber \\
     &   &  +\,\, 2\,102\,099x^3 + 1\,877\,667x^2 + 3\,904\,161x + 1\,164\,691)^2,
           \nonumber \\
R_d(x) & = & 2^6\,3^{15} (5x^6 - 6x^5 + 111x^4 + 64x^3 + 795x^2 + 1254x + 5477). 
           \nonumber
\end{eqnarray}

By the way, a naive approach mentioned in Sect.~\ref{sec:computation},
namely, taking polynomials $A$ and $B$ of degrees 10 and 15 respectively with 
indeterminate coefficients and equating to zero the coefficients of degrees 
from 7 to 30 of $A^3-B^2$, would, this time, lead us to a system of polynomial 
equations of degree $6\,198\,727\,824$. It took 40 years (from 1965 to 2005)
to compute all the four DZ-pairs of this example, but the fact that there are
exactly four non-equivalent solutions and that two of them are defined over $\Q$
while the other two are defined over an imaginary quadratic field, can be
immediately seen from the picture without any computation.


\section{Some sporadic examples of Beukers and \\ Stewart \cite{BeuSte-10}}\label{sec:BS}

All the polynomials in this section which correspond to the {\em asymmetric}\/ trees
are taken from the above-cited article \cite{BeuSte-10}. The normalization
sometimes is changed. The goal of this section is to show the 
{\em combinatorial reasons}\/ of appearance of these sporadic examples.

\subsection{Passport $(7^3,3^7)$}\label{sec:split1}

\begin{figure}[htbp]
\begin{center}
\epsfig{file=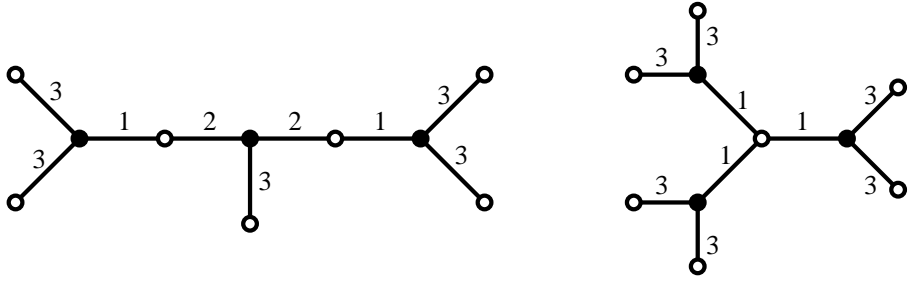,width=12cm}
\caption{Two trees corresponding to the passport $(7^3,3^7)$; one of them 
is symmetric, the other one is not. Therefore, both are defined over $\Q$.}
\label{fig:split-7-3}
\end{center}
\end{figure}

The passport shows that we are treating here the problem of the minimum degree of
the difference $A^7-B^3$ where $\deg A=3$, $\deg B=7$.
The combinatorial orbit consists of two trees, see Fig.~\ref{fig:split-7-3}. One of them 
is symmetric, the other one is not; therefore, both are defined over $\Q$.

\ni
The triple corresponding to the asymmetric tree is as follows:
\begin{eqnarray}
P & = & (x^3+18x+18)^7, \nonumber \\
Q & = & (x^7+42x^5+42x^4+504x^3+1008x^2+1512x+3024)^3 \nonumber \\
R & = & 2^4\,3^3\, (77x^{12}+5922x^{10} + 6237x^9 + 172\,368x^8 + 366\,606x^7 + 2\,451\,330x^6
	    \nonumber \\
  &   & +\,\, 7\,314\,300x^5 + 19\,105\,632x^4 + 53\,867\,268x^3 + 82\,260\,360x^2 \nonumber \\
  &   & +\,\, 86\,097\,816x + 62\,594\,856). \nonumber
\end{eqnarray}

\ni
The triple corresponding to the symmetric tree may be computed as follows:
\begin{enumerate}
\item   Compute the polynomials corresponding to a branch of this three-branch tree, 
		that is, to a tree of the series $A$ (see Sect.~\ref{sec:A}) with the
		parameters $s=3$, $t=1$, $k=2$.
\item	Make the change of variables $x\to 1-x$ in order to put the white vertex of
		degree~1 to the point~$x=0$; thus, the polynomial $P(x)$, instead of being~$x^7$, 
		becomes $(1-x)^7$; it is convenient to change its sign and to get $(x-1)^7$.
\item   Insert $x^3$ instead of $x$.
\end{enumerate}
By pure convenience we add to the above operations one more: instead of taking 
$P(x)=(x-1)^7$ we take $P(x)=(x-3)^7$. This permits us to avoid fractional coefficients. 
The resulting triple is
\begin{eqnarray}
P & = & (x^3-3)^7, \nonumber \\ 
Q & = & x^3(x^6-7x^3+14)^3, \nonumber \\
R & = & -14x^{12}+189x^9-987x^6+2359x^3-2187. \nonumber
\end{eqnarray}

\subsection{Passport $(8^3,3^8)$}\label{sec:split2}

\begin{figure}[htbp]
\begin{center}
\epsfig{file=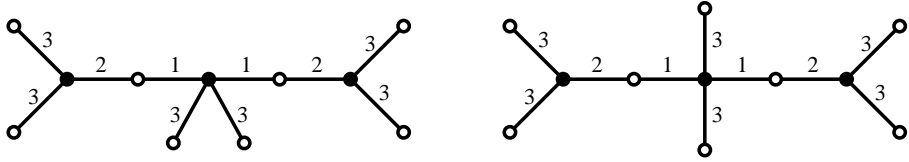,width=12cm}
\caption{Two trees corresponding to the passport $(8^3,3^8)$.} 
\label{fig:split-8-3}
\end{center}
\end{figure}

The passport corresponds to the problem of the minimum degree of the difference
$A^8-B^3$ where $\deg A=3$, $\deg B=8$.
The combinatorial orbit consists of two trees, see Fig.~\ref{fig:split-8-3}. One of them 
is symmetric, the other one is not; therefore, both are defined over $\Q$.

\ni
The triple corresponding to the asymmetric tree looks as follows:
\begin{eqnarray}
P & = & (x^3+27x+81)^8, \nonumber \\
Q & = & (x^8+72x^6+216x^5+1620x^4+9720x^3+24300x^2+87480)^3, \nonumber \\
R & = & -3^{10}\,(52x^{14}+6942x^{12}+21\,816x^{11}+366\,444x^{10}+2\,319\,840x^9 \nonumber \\
  &   & +\,\,13\,129\,047x^8+90\,716\,760x^7+406\,062\,720x^6+1\,812\,830\,544x^5 \nonumber \\
  &   & +\,\,7\,862\,190\,642x^4+23\,694\,237\,936x^3+67\,352\,942\,772x^2 \nonumber \\
  &   & +\,\,173\,534\,618\,376x+204\,401\,597\,391). \nonumber
\end{eqnarray}

\ni
The triple corresponding to the symmetric tree may be computed as follows:
\begin{enumerate}
\item   Compute the polynomials corresponding to the series $E_4$ (see Sect.~\ref{sec:E-even})
		with $s=1$, $t=2$, $k=1$, $l=2$.
\item	Make the change of variables $x\to x+1$ in order to move the (left) black vertex
		of degree~4 from $-1$ to 0.
\item   Insert $x^2$ instead of $x$.
\end{enumerate}
We omit the resulting polynomials.

\subsection{Passport $(10^3,3^{10})$}\label{sec:split3}

\begin{figure}[htbp]
\begin{center}
\epsfig{file=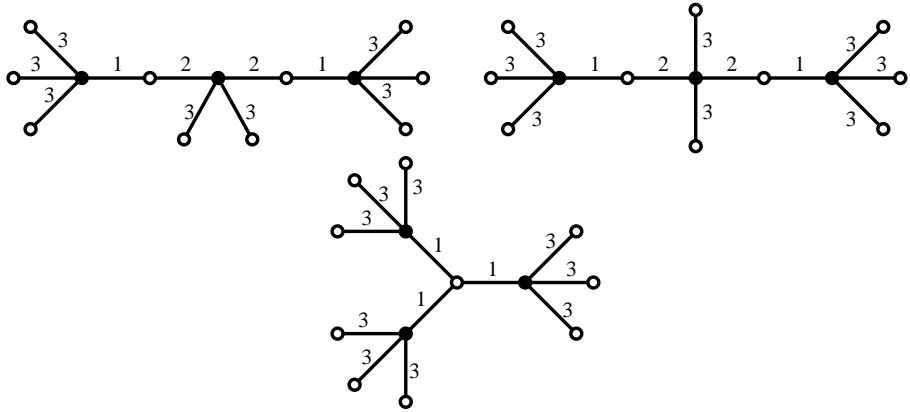,width=12cm}
\caption{Three trees corresponding to the passport $(10^3,3^{10})$.} 
\label{fig:split-10-3}
\end{center}
\end{figure}

This time we deal with the problem $\min\deg(A^{10}-B^3)$, $\deg A=3$, $\deg B=10$.
The combinatorial orbit corresponding to this passport contains three trees,
see Fig.~\ref{fig:split-10-3}. These trees have three different symmetry
types, hence all of them are defined over $\Q$.

\ni
The polynomials for the asymmetric tree look as follows:
{\small
\begin{eqnarray}
P & = & (x^3+54x+162)^{10}, \nonumber \\
Q & = & (x^{10}+180x^8+540x^7+11\,340x^6+68\,040x^5+374\,220x^4 \nonumber \\
  &   & +\,\,2\,449\,440x^3+8\,573\,040x^2+22\,044\,960x+57\,316\,896)^3, \nonumber \\
R & = & -2^4\,3^{11}\,(595x^{18}+201\,960x^{16}+629\,748x^{15}+28\,669\,140x^{14} \nonumber \\
  &   & +\,\,179\,596\,440x^{13}+2\,460\,946\,860x^{12}+20\,601\,540\,000x^{11} \nonumber \\
  &   & +\,\,158\,558\,654\,736x^{10}+1\,257\,674\,415\,840x^9+7\,823\,104\,403\,040x^8 \nonumber \\
  &   & +\,\,46\,607\,404\,043\,520x^7+253\,091\,029\,021\,200x^6+
        1\,120\,772\,437\,834\,752x^5 \nonumber \\
  &   & +\,\,4\,520\,664\,857\,839\,680x^4+15\,435\,507\,254\,345\,280x^3 \nonumber \\
  &   & +\,\,37\,331\,470\,988\,020\,800x^2+62\,014\,139\,393\,904\,000x \nonumber \\
  &   & +\,\,62\,042\,237\,538\,382\,656). \nonumber
\end{eqnarray}
}

\vspace{-2mm}

The polynomials for the tree with the symmetry of order 2 is computed in the same way
as in Sect.~\ref{sec:split2}. The parameters of the tree of the type $E_4$ are $s=2$, 
$t=1$, $k=1$, $l=3$; then we must replace $x$ with $x+1$, and insert $x^2$ instead of $x$.

The polynomials for the tree with the symmetry of order 3 is computed in the same way
as in Sect.~\ref{sec:split1}. The parameters of the tree of the type $A$ are $s=3$, 
$t=1$, $k=3$; then we must replace $x$ with $1-x$, and insert $x^3$ instead of~$x$.


\subsection{Passport $(9^5,5^9)$}\label{sec:why}

\begin{figure}[htbp]
\begin{center}
\epsfig{file=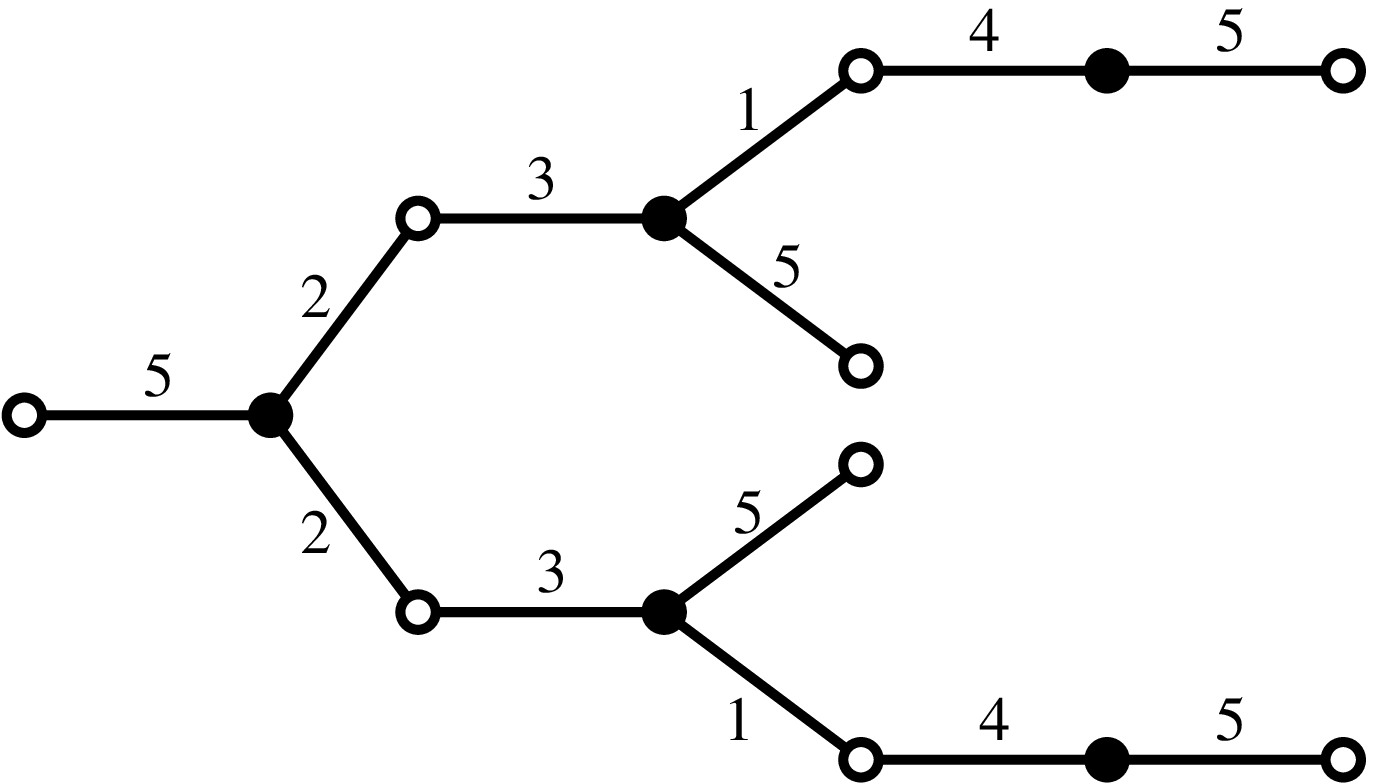,width=7.2cm}
\caption{This tree, corresponding to the passport $(9^5,5^9)$,
is defined over~$\Q$. All known combinatorial invariants of Galois
action fail to explain this phenomenon.} 
\label{fig:why}
\end{center}
\end{figure}

We finish this section with an example which shows that the combinatorial methods,
while being very powerful, are, however, not all-powerful. There are several trees
with the passport $(9^5,5^9)$, and one of them, shown in Fig.~\ref{fig:why}, is
defined over $\Q$ without any apparent reason. All known combinatorial and 
group-theoretic Galois invariants fail to explain this phenomenon. All we can say
is that the corresponding system has rational solutions ``by chance''.

\msk

\ni
The polynomials $P$ and $Q$ for the tree of Fig.~\ref{fig:why} are as follows:
\begin{eqnarray}
P & = & (x^5+50x^3+500x+500)^9, \nonumber \\
Q & = & (x^9+90x^7+2700x^5+900x^4+30\,000x^3+36\,000x^2 \nonumber \\
  &   & +\,\,90\,000x+180\,000)^5. \nonumber
\end{eqnarray}
The polynomial $R$ here is of degree 32, and it is too cumbersome, so we do not
write it explicitly.

\section{Yet more examples}\label{sec:other}

\subsection{An infinite series of splitting combinatorial orbits}

We have already seen two examples (see Sect.~\ref{sec:split1} and \ref{sec:split2})
of combinatorial orbits of size~2 which, instead of being defined over a quadratic
field, split in two orbits defined over $\Q$ because the trees in question have
different orders of symmetry. Here we present an infinite series of such examples.

\begin{figure}[htbp]
\begin{center}
\epsfig{file=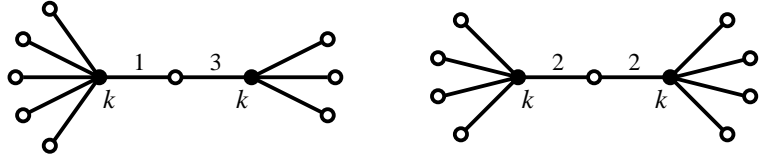,width=10cm}
\end{center}
\caption{Two trees with the passport $(k^2,4^11^{2k-4})$. One of them is
symmetric, the other one is not.}
\label{fig:sym-and-not}
\end{figure}

The trees in question have the passport $(k^2,4^11^{2k-4})$ for $k\ge 3$, see
Fig.~\ref{fig:sym-and-not}. Belyi function for the symmetric tree looks as follows:
\begin{eqnarray}\label{sym}
f_1(x) \eq \frac{(-1)^{k+1}}{k^k}\cdot\frac{(x^2-k)^k}{x^2-1}. \nonumber
\end{eqnarray}
Belyi function for the asymmetric tree looks as follows:
{\small
$$
f_2(x) \eq \frac{(-1)^k}{(6k)^{k-1}(k-2)^{k-2}(2k-1)^{2k-1}}\cdot
\frac{(x^2-6k(2k-1)x-6k(k-2)(2k-1)^2)^k}{x^2+6k(k-2)x+6k(k-2)^2(2k-1)}\,.
$$
}

In both cases, the white vertex of degree 4 lies at $x=0$.
The expressions for Belyi functions give us the polynomials $P$ and $R$. 

In order to prove the correctness of the above expressions we need to verify
two things: for both $f_1$ and $f_2$, we have (a) $f(0)=1$; (b) first three 
derivatives of $f(x)$ at $x=0$ vanish. 

We leave the proof to the reader.


\subsection{Trees with a relaxed minimum degree condition}

Let us return to the problem of the minimum degree of the difference
$A^3-B^2$, the question from which this whole line of research started
(see \cite{BCHS-65}. We have seen that when $\deg A=2k$, $\deg B=3k$,
we have $\min\deg(A^3-B^2)=k+1$. For $k\ge 6$, the computation becomes
exceedingly difficult, and there is practically no hope to find solutions
defined over $\Q$. However, if we are not so demanding and accept a solution
with the degree of $A^3-B^2$ slightly greater than $k+1$, then sometimes
we can find a needed solution.

\begin{figure}[htbp]
\begin{center}
\epsfig{file=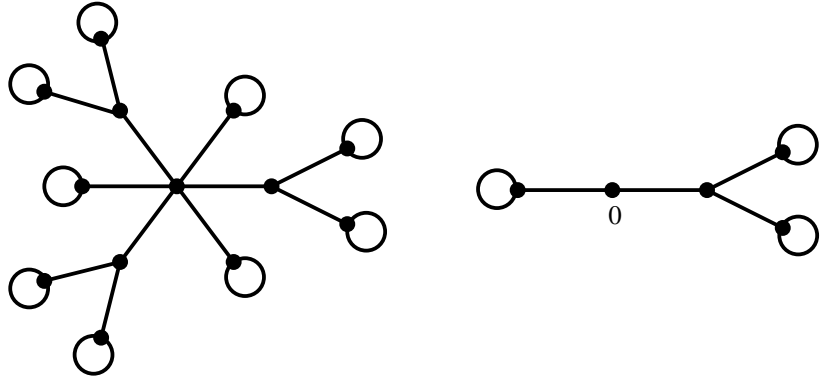,width=10.8cm}
\caption{The map on the left represents two polynomials $A$ and $B$, of degrees
$2k=14$ and $3k=21$ respectively, such that $\deg\,(A^3-B^2)=9$. Thus, 
the degree of the difference does not attain its minimum value $k+1=8$, 
but in return both $A$ and $B$ are defined over $\Q$.} 
\label{fig:deg9-instead-8}
\end{center}
\end{figure}

\begin{example}
Let us take a polynomial $A$ with one double root, so that $A^3$ would have
one root of multiplicity 6 and all the other roots of multiplicity 3.
The corresponding tree would have one vertex less and therefore one face more.

The tree in Fig.~\ref{fig:deg9-instead-8} corresponds to $k=7$. It is the ``cube'' 
of the tree $S$, see Sect.~\ref{sec:tree-S}. Therefore, all we have to do is to 
insert $x^3$ instead of $x$ in the formulas of Section~\ref{sec:tree-S}.
\begin{eqnarray}
P & = & x^6\,(x^{12}+24x^9+176x^6-2816)^3, \nonumber \\
Q & = & (x^{21}+36x^{18}+480x^{15}+2304x^{12}-3840x^9 -55\,296x^6 \nonumber \\
  &   & -\,\,14\,336x^3+221\,184)^2, \nonumber \\
R & = & 2^{22}\cdot 3^3\,(x^9+17x^6+56x^3-432). \nonumber
\end{eqnarray}
\end{example}

\begin{figure}[htbp]
\begin{center}
\epsfig{file=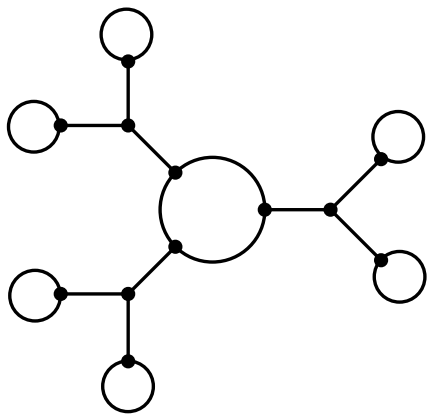,width=5.6cm}
\caption{This map represents two polynomials $A$ and $B$, of degrees
$2k=12$ and $3k=18$ respectively, such that $\deg\,(A^3-B^2)=9$. Thus, 
the degree of the difference does not attain its minimum value $k+1=7$, 
but in return both $A$ and $B$ are defined over $\Q$.} 
\label{fig:deg-relaxed}
\end{center}
\end{figure}

\begin{example}
When all the roots of $A$ and $B$ are distinct, the polynomial $R$ has
$k+1$ distinct roots. Let us accept $R$ with a multiple root (thus, its degree
will be greater that $k+1$). The tree in Fig.~\ref{fig:deg-relaxed} gives
such and example. It corresponds to $k=6$, and $\deg R=9$.
The polynomials for this tree look as follows:
\begin{eqnarray}
P & = & (x^3+3)^3\,(x^9+9x^6+27x^3+3)^3, \nonumber \\
Q & = & (x^{18}+18x^{15}+135x^{12}+504x^9+891x^6+486x^3-27)^2, \nonumber \\
R & = & 1728x^3\,(x^6+9x^3+27). \nonumber
\end{eqnarray}
\end{example}

\paragraph{Acknowledgements.} Fedor Pakovich was partially supported by ISF grants 
No.~639/09 and 779/13. He is also grateful to the Max Plank Institute for Mathematics 
for the hospitality and support. Alexander Zvonkin was partially supported by the 
Research grant {\sc Graal} ANR-14-CE25-0014.

\vspace{1cm}

\ni
{\sc Fedor Pakovich:} Department of Mathematics, Ben-Gurion University
of the Negev, P.O.B. 653, Beer Sheva, Israel; e-mail: {\tt pakovich@math.bgu.ac.il}

\msk

\ni
{\sc Alexander K. Zvonkin:} LaBRI, UMR 5800, Universit\'e de Bordeaux, 33400 Talence, 
France; e-mail: {\tt zvonkin@labri.fr}, and The Chebyshev Mathematical Laboratory at the
Saint-Petersburg State University, 29B, 14th Line, Vasilyevsky Island, Saint-Petersburg 
199178, Russia.

\end{document}